\documentclass[12pt,a4paper,french,draft]{smfart}
\usepackage{amssymb,amsmath,amscd,stmaryrd,vmargin}
\usepackage[francais]{babel}
\usepackage[all]{xy}
\usepackage[T1]{fontenc}
\usepackage[utf8]{inputenc}
\theoremstyle{plain}
\newtheorem{theo}{Théorème}[section]
\newtheorem{prop}[theo]{Proposition}
\newtheorem{cor}[theo]{Corollaire}
\newtheorem{lem}[theo]{Lemme}
\newtheorem{defi}[theo]{Définition}

\def \plnc {\phi_{\ell,n}^\chi}
\def \pln {\phi_{\ell,n}}
\def \lnis {\LM_n^\text{\rm iso}}
\def \cfrc {C_{F,r}^\chi}
\def \hsrc {H_{S,r}^\chi}

\def \demdu#1 { {\sl Proof #1.} }
\def \Hom {\text{\rm Hom}}
\def \Gal {\text{\rm Gal}}

\def \limpro{\lim\limits_{\leftarrow} }

\def \qed{\text {$\quad \square$}}

\def \U {\overline {U}}

\def \F#1 {\overline {F^\times_{#1}}}

\def\CM{{\mathbb{C}}}

\def\FM{{\mathbb{F}}}

\def\LM{{\mathbb{L}}}

\def\NM{{\mathbb{N}}}

\def\QM{{\mathbb{Q}}}

\def\ZM{{\mathbb{Z}}}

\def\al{\alpha}

\def\ga{\gamma}

\def\de{\delta}
\def\De{\Delta}
\def\ep{\varepsilon}

\def\la{\lambda}
\def\La{\Lambda}

\def\si{\sigma}

\def\ze{\zeta}

\def\EC{{\mathcal{E}}}

\def\GC{{\mathcal{G}}}

\def\LC{{\mathcal{L}}}

\def\RC{{\mathcal{R}}}

\def\UC{{\mathcal{U}}}
\def\VC{{\mathcal{V}}}

\def\xba{{\overline{x}}}

\def\Cba{{\overline{C}}}

\def\chit{{\widetilde{\chi}}}

\def\Gal{\mathop{\mathrm{Gal}}\nolimits}
\def\Hom{\mathop{\mathrm{Hom}}\nolimits}

\def\Ker{\mathop{\mathrm{Ker}}\nolimits}
\def\ker{\Ker}

\providecommand{\bysame}{\leavevmode\hbox to3em{\hrulefill}\thinspace}
\providecommand{\MR}{\relax\ifhmode\unskip\space\fi MR } 
\providecommand{\MRhref}[2]{}
\providecommand{\href}[2]{#2}
\providecommand{\bysame}{\leavevmode ---\ }
\providecommand{\og}{``}
\providecommand{\fg}{''}

\providecommand{\smfandname}{et}

 \def\Dbar{\leavevmode\lower.6ex\hbox to
0pt{\hskip-.23ex \accent"16\hss}D}
  \def\cfac#1{\ifmmode\setbox7\hbox{$\accent"5E#1$}\else
  \setbox7\hbox{\accent"5E#1}\penalty 10000\relax\fi\raise 1\ht7
  \hbox{\lower1.15ex\hbox to 1\wd7{\hss\accent"13\hss}}\penalty 10000
\hskip-1\wd7\penalty 10000\box7}
\def\cftil#1{\ifmmode\setbox7\hbox{$\accent"5E#1$}\else
  \setbox7\hbox{\accent"5E#1}\penalty 10000\relax\fi\raise 1\ht7
  \hbox{\lower1.15ex\hbox to 1\wd7{\hss\accent"7E\hss}}\penalty 10000
\hskip-1\wd7\penalty 10000\box7}

\title[]{Indices isotypiques des éléments cyclotomiques.}    \author{Tatiana Beliaeva}
\address{IRMA, UMR 7501 de l'Université Louis Pasteur et du CNRS\\
7 rue René-Descartes, 67084 Strasbourg Cedex, France
}
\email{beliaeva@math.u-strasbg.fr}
\author{Jean-Robert Belliard}
\address{ Laboratoire de Mathématiques,
UMR 6623 de l'Université de Franche-Comté et du CNRS\\
16 route de Gray, 25030 Besan\c{c}on cedex, France}
\email{jean-robert.belliard@univ-fcomte.fr}

\date{\today}
\begin{document}

\frontmatter
\begin{abstract}
Given $F$ a real abelian field, $p$ an odd prime and $\chi$ any Dirichlet character of $F$, we give a method for computing the $\chi$-index $\displaystyle \left (H^1(G_S,\ZM_p(r))^\chi: C^F(r)^\chi\right)$ where the Tate twist $r$ is an odd integer $r\geq 3$, the group $C^F(r)$ is the group of higher circular units, $G_S$ is the Galois group over $F$ of the maximal $S$ ramified algebraic extension of $F$, and $S$ is the set of places of $F$ dividing $p$. This $\chi$-index can now be computed in terms only of elementary arithmetic of finite fields $\FM_\ell$. Our work generalizes previous results by Kurihara who used the assumption that the order of $\chi$ divides $p-1$.
\end{abstract}

\maketitle
\tableofcontents

\mainmatter

\section*{Introduction}
Soit $F$ un corps de nombres abélien totalement réel, d'anneau d'entiers $O_F$, de groupe de Galois $G=\Gal(F/\QM)$ et $p\neq 2$ un nombre premier rationnel fixé dans toute la suite. On note $S=S(F)$ l'ensemble des places de $F$ divisant $p$, on note $\Omega_S$ l'extension algébrique maximale, non ramifiée hors de $S$, de $F$ et $G_S=\Gal(\Omega_S/F)$. Les objets mathématiques centraux de cet article sont les groupes de cohomologie galoisienne $H^i(G_S,\ZM_p(r))=\limpro H^i(G_S,\ZM/p^n(r))$, pour $r\in \NM$, $r\geq 2$ et pour $i=1$ ou $2$, vus comme $\ZM_p[G]$-modules. Ces groupes de cohomologie galoisienne sont connus pour être isomorphes aux groupes de cohomologie étale $H^i(\mathrm{spec}\, O_F[1/p],\ZM_p(r))$. Parmi les motivations à leur étude il y a leur interprétation en $K$-théorie. Sous la conjecture de Quillen-Lichtenbaum ces groupes sont aussi isomorphes  pour $r\geq 2$ aux $p$-adifiés $K_{2r-i} (O_F)\otimes \ZM_p$ des groupes de $K$-théorie de Quillen (voir \cite{Kosur75} pour la construction des groupes de $K$-théorie supérieurs). Cependant notre travail se situe entièrement en cohomologie galoisienne $p$-adique et est donc logiquement indépendant de ce contexte qui était conjectural au début de la rédaction de ce travail. Entre temps les travaux de Voevodsky, complétés par Weibel, ont aboutit à une démonstration de cette conjecture de Quillen-Lichtenbaum (voir \cite{V10}). 

En arithmétique ces groupes de cohomologie $p$-adique s'interprètent comme des analogues supérieurs des groupes des $p$-unités (pour $i=1$ et $r\not\in \{0,1\}$) et des $p$-groupes de $S$-classes d'idéaux (pour $i=2$ et $r\geq 0$). Ces analogies sont précisées par une suite exacte due à Kolster, {N}guy{\cftil{e}}n-{Q}uang-{D}{\cftil{o}} et Fleckinger pour $i=1$ (voir theorem 3.2 de \cite{KNF} p.689) et une suite exacte due à Kurihara pour $i=2$ (voir lemma 4.3 de \cite{Kur99} p. 269). Ils sont aussi particulièrement intéressants à étudier en liaison avec les valeurs spéciales des fonctions $L$ en $s=1-r$ : ils sont au premier plan dans les démonstrations des conjectures à la Bloch-Kato pour les valeurs spéciales des fonctions $L$ motiviques dans le cas particulier des fonctions $L$ de Dirichlet (voir entre autres \cite{BK90,KNF,BenN02,HK03,BG03}).

Un des points culminants dans l'article \cite{KNF} est une formule d'indice pour $F$ totalement réel abélien et $r>1$ impair, conséquence du theorem 5.4 de \cite{KNF}~:
$$\# \left ( \frac {H^1(G_S,\ZM_p(r))} {C^F(r)} \right ) =\# H^2(G_S,\ZM_p(r))$$
Ici $C^F(r)$ (voir la définition \ref{CFr}) est l'analogue supérieur des unités circulaires, et cette formule d'indice est donc l'analogue supérieur de la formule d'indice de Sinnott : une version algébrique de la formule analytique du nombre de classes de Dedekind.
Cette formule d'indice passe aux $\chi$-composantes lorsque $p\nmid [F:\QM]$, pour tout caractère $\chi$ sur $G=\Gal(F/\QM)$ (c'est aussi une conséquence du theorem 5.4 de \cite{KNF}). Lorsque $p\mid [F:\QM]$ l'indice \(\displaystyle \left (H^1(G_S,\ZM_p(r))^\chi: C^F(r)^\chi\right) \) et l'ordre \( \# H^2(G_S,\ZM_p(r))^\chi\) sont très probablement reliés, mais il y a a priori une déviation difficile à expliciter. Le principal résultat de cet article est de présenter une méthode \underline{explicite} et \underline{élémentaire} pour calculer les $\chi$-indices $\displaystyle \left (H^1(G_S,\ZM_p(r))^\chi: C^F(r)^\chi\right)$ pour \underline{tous} les caractères de Dirichlet $\chi$ (voir le théorème \ref{main}, son corollaire \ref{pgcd} et la formulation explicite du  théorème \ref{formexpire}).

Dans l'article \cite {Kur99} Kurihara introduit une méthode remarquable pour calculer ces $\chi$-indices en liaison avec l'ordre du $H^2(G_S,\ZM_p(r))$, qu'il considère comme le $p$-adifié du groupe de Tate-Shafarevich associé au motif $h^0(\QM(\zeta_f))$. La méthode de Kurihara ramène le calcul de $\chi$-indice à de l'arithmétique dans des corps finis $\FM_\ell$, mais n'est complètement décrite dans l'article \cite{Kur99} que pour les caractères de Dirichlet $\chi$ dont l'ordre divise $p-1$. L'apport principal du présent article est de se dispenser de cette hypothèse restrictive. Pour ce faire on est parvenu à remplacer les considérations sur l'ordre de certains éléments spéciaux dans \cite{Kur99} par des isomorphismes canoniques et généraux de $\chi$-composantes
(voir le théorème \ref{inj}).

Cet article est composé comme suit. Dans la section \ref{mordes} on définit les notations $H^1(G_S,\ZM_p(r))$ en justifiant les restrictions $r\geq 3$ impairs et $F$ totalement réel. Puis on rappelle la construction d'un morphisme de descente classique (voir par exemple \cite{KNF}) $\U'_\infty(r-1)_{G_\infty}\overset \al \longrightarrow H^1(G_S,\ZM_p (r))$ où $\U'_\infty$ désigne la limite projective pour les normes des $p$-adifiés des $p$-unités des $F_n=F(\mu_{p^n})$ le long des étages finis de la tour $F_\infty=\bigcup_n F_n$ et où $G_\infty =\Gal(F_\infty/F)$. Il s'agit surtout, pour l'usage ultérieur fait dans cet article, de préciser explicitement les images de $p$-unités données.

Dans la section \ref{eltwis} on définit, en utilisant le morphisme $\al$ et suivant \cite{KNF} et \cite{Kur99}, les éléments cyclotomiques en cohomologie galoisienne $c^F_b(r)$ qui engendrent le groupe $C^F(r)$ et coïncident avec le $p$-adifiés des éléments cyclotomiques de Soulé (\cite{Sou80} ou aussi \cite{Sou87}) sous la conjecture de Quillen-Lichtenbaum. On étudie aussi les relations de normes reliant ces générateurs. Suivant le même esprit qu'en section \ref{mordes} on précise explicitement les images de ces $c^F_b(r)$ par restriction et réduction modulo $p^n$ dans $H^1(G_S(F_n),\ZM/p^n (r))$.

Dans la section \ref{carut} on calcule le caractère commun (qui est le caractère régulier) des $G=\Gal (F/\QM)$-représentations linéaires obtenues à partir de $C^F(r)\otimes \QM_p$, $ H^1(G_S,\QM_p (r))$ et $\U'_\infty(r-1)_{G_\infty}\otimes \QM_p$. Pour ce faire on utilise l'injectivité et la finitude du conoyau de $\al$ démontrées dans \cite{KNF}; et la théorie de Coleman (\cite{Col79}) pour calculer le caractère des unités semi-locales et en déduire celui des $p$-unités.

Dans la section \ref{chipart} on regroupe les lemmes algébriques concernant les $\chi$-parties qui servent à notre étude. Pour des descriptions plus complètes des propriétés fonctorielles des $\chi$-parties et $\chi$-quotients, le lecteur est invité à consulter \cite{So90}, ou encore \cite{T99}. On indique aussi quel $\chi$-indice précis est calculé dans la suite, parmi la multitude de nuances possibles lorsque $p\mid [F:\QM]$. Tant que faire se peut, on justifie le choix de $\chi$-indice fait dans cet article.

Dans la section \ref{redl} on reprend la construction par Kurihara du morphisme $\phi_{\ell,n}$ de réduction modulo $\ell$ qui est au c\oe ur de cette méthode de calcul des $\chi$-indices. On détaille soigneusement une interprétation kummerienne de ce morphisme avec le diagramme-clé (\ref{kumres}); et l'on démontre le théorème d'isomorphisme \ref{inj} qui est le point essentiel dans notre amélioration de l'approche de Kurihara.

Dans la section \ref{thmind} on utilise le théorème \ref{inj} pour démontrer la formule d'indice du théorème \ref{main} qui est le résultat principal de cet article. Par rapport à la démarche de \cite{Kur99}, pour démontrer ce théorème il faut essentiellement dans cette section circonvenir les difficultés techniques  supplémentaires occasionnés par l'éventualité $p\mid [F:\QM]$  (cas non semi-simple). C'est ici qu'interviennent à la fois les propriétés fonctorielles des $\chi$-parties rappelées en section \ref{chipart} et le gain conceptuel obtenu avec l'isomorphisme du théorème \ref{inj}.

Dans la section \ref{calcexp} on précise l'énoncé du théorème \ref{main} en donnant une formule explicite pour les images d'éléments cyclotomiques contre le morphisme $\phi_{\ell,n}$ de réduction modulo $\ell$ et en fixant des générateurs concrets des $\chi$-parties considérées.

Dans la section \ref{numerique} on donne quelques exemples num\'eriques dans le cas ou $F$ est un corps de degr\'e $p$ et $\chi$ est un caract\`ere d'ordre $p$.

\noindent{\it Remerciements :} Les deux auteurs remercient Thong  {N}guy{\cftil{e}}n-{Q}uang-{D}{\cftil{o}} qui nous a suggéré ce sujet comme une collaboration possible. La section \ref{numerique} a été ajoutée à la demande du referee anonyme, nous le remercions aussi pour sa lecture attentive et ses conseils constructifs. 
Ce travail doit énormément à Kurihara. D'abord l'article \cite{Kur99} nous a inspiré. Ensuite nous avons profité de ses conseils lors de la mise en forme finale. Nous lui en sommes particulièrement reconnaissant.

\section{Morphisme de descente explicite}\label{mordes}



On fixe un plongement de $\overline {\QM}$ dans $\CM$. Ce plongement définit un système de racines de l'unité compatibles en posant $\displaystyle \ze_n:=e^{2i\pi/n}$, pour $n\in \NM$, $n\geq 1$. On rappelle que $F$ est un corps de nombres abélien totalement réel et $p\neq 2$ est un nombre premier rationnel fixé dans toute la suite.
 
Soit $S=S(F)$ l'ensemble de places de $F$ au-dessus de $p$. On note $\Omega_S$ l'extension algébrique maximale de $F$ non ramifiée en dehors de $S$  et $G_S=\Gal(\Omega_S/F)$ son groupe de Galois. On note aussi $\GC_S=\Gal(\Omega_S/\QM)$.
Soit $F_n=F(\ze_{p^n})\subset \Omega_S$. On note $G_S(F_n)$ le groupe $\Gal(\Omega_S/F_n)$, on note $G_n=\Gal(F_n/F)$ et $G_\infty=\Gal(F_\infty/F)\cong \limpro G_n$, où $F_\infty = \bigcup_n F_n$.

 La donnée du système compatible en norme $\ze_{p^n}$ définit aussi un générateur $t(1)$ du module de Tate $\ZM_p(1):=\limpro \mu_{p^n}$, où $\mu_m$ désigne le groupe des racines de l'unité d'ordre divisant $m$ dans $\overline{\QM}$. Pour tout entier $r\geq 1$ on note $t(r)\in \ZM_p(r)$ la puissance tensorielle $r^{\text{ième}}$ de $t(1)$ qui est un générateur de $\ZM_p(r):=\ZM_p(1)^{\otimes^r}$ et pour tout entier $n\geq 1$ on note $t(r)_n$ la $n^{\text{ième}}$ projection de $t(r)$ qui engendre $\ZM /p^n (r)$. Si $\kappa \colon \GC_S \twoheadrightarrow \ZM_p^\times$ désigne le "caractère" cyclotomique, alors l'action de $\GC_S$ sur $\ZM_p(r)$ est donnée par $\ga ( x) = \kappa(\ga)^r x$. Ces conventions permettent de définir pour tout entier $r\geq 1$ et tout $\ZM_p[\GC_S]$-module $M$ son $r^{\text{ième}}$-tordu à la Tate, noté $M(r)$, par $M(r):= M\otimes_{\ZM_p} \ZM_p(r)$, le groupe $\GC_S$ agissant diagonalement. Il est parfois plus commode de voir $M(r)$ comme le $p$-groupe $M$ lui-même muni de l'action de $\GC_S$
tordue : $\ga *_r x := \kappa(\ga)^r \ga (x)$. Le principal objet étudié dans cet article est le groupe de cohomologie galoisienne~: $$ H^1(G_S,\ZM_p(r)) := {\limpro}_{n} H^1(G_S, \ZM/p^n (r)).$$
Suivant Kurihara (\cite{Kur99}) on va se limiter aux twists $r$ impairs, à un corps $F$ totalement réel et supposer $r > 1$. Cette approche est basée sur l'existence d'éléments spéciaux engendrant un sous-$\ZM_p[G]$-module libre de rang $1$ et d'indice fini dans $H^1(G_S,\ZM_p(r))$ (voir théorème \ref{freespe}). D'après les calculs de rang p.~223 de \cite{Kosur} la seule possibilité pour que  $H^1(G_S,\ZM_p(r))$ soit de $\ZM_p$-rang égal à l'ordre de $G=\Gal (F/\QM)$ est de prendre le corps $F$ totalement réel et le twist $r>1$ impair.
Dans le cas du twist $r=1$ on a $H^1(G_S,\ZM_p(r)) \cong \U'_F:=U'_F\otimes \ZM_p$, la pro-$p$-complétion du groupe des $S$-unités de $F$. Le twist $r=1$ fait l'objet d'un nombre considérable de travaux, est connu pour être techniquement plus difficile et ne sera pas abordé dans cet article.
 L'un des avantages plus accessoires de ces restrictions en généralité est le lemme suivant~:
\begin{lem}\label{torfree} On suppose $r>1$ impair, $F$ totalement réel et $p\neq 2$; alors le $\ZM_p$-module $H^1(G_S,\ZM_p(r))$ est sans torsion.
\end{lem}
\begin{proof} D'après le lemme 2.2 de \cite{KNF} la torsion de $H^1(G_S,\ZM_p(r))$ s'identifie à $\QM_p/\ZM_p(r)^{G_\infty}$. Comme $F$ est totalement réel la conjugaison complexe $\tau$ est un élément bien défini de $G_\infty$ qui vérifie $\kappa(\tau)^r = (-1)^r=-1$ puisque $r$ est impair. Cet élément agit donc simultanément par multiplication par $1$ et $-1$ sur $\QM_p/\ZM_p(r)^{G_\infty}$ ce qui démontre que la torsion de $H^1(G_S,\ZM_p(r))$ est annulé par $2$.
\end{proof}
Pour finir cette section on redonne la construction et une formule explicite pour un morphisme de descente
classique (voir \cite{KNF} par exemple)~:
$$\U'_\infty(r-1)_{G_\infty}\overset \al \longrightarrow H^1(G_S,\ZM_p (r)),$$
où $\U'_\infty$ désigne la limite projective relativement aux applications de norme des pro-$p$-complétions des $S$-unités $\U'_n=\U'_{F_n}$ le long de l'extension $F_\infty/F$.

Soit $u_\infty=(u_n)_{n\in\NM}$ une suite cohérente en norme de $\U'_\infty$ avec $u_n\in  F_n$. Alors par le cocycle de Kummer $\U'_n/p^n \otimes \ZM/p^n (r-1) \hookrightarrow H^1(G_S(F_n),\ZM/p^{n} (r))$, chaque $u_n \otimes t(r-1)_{n}$ définit un élément $x_n\in H^1(G_S(F_n),\ZM/p^{n} (r))$. Concrètement on a la formule $$x_n(g)=\frac {g(\sqrt[p^n]{u_n})} {\sqrt[p^n]{u_n}} \otimes t(r-1)_n \in \mu_{p^n}(r-1) \cong \ZM/p^n (r).$$ Pour tout $m$ on pose
$$
\al_{m} (u_\infty\otimes t(r-1))= cor_{F(\ze_{p^{m}})/F} (x_m) \in H^1(G_S,\ZM/p^{m}(r)).
$$
La cohérence en norme des $u_n$ impose alors que les $\al_m(u_\infty)$ forment une suite cohérente pour les applications $H^1(G_S,\ZM/p^{m+k}(r)) \longrightarrow H^1(G_S,\ZM/p^{m}(r))$, qui définit donc un élément de $ H^1(G_S,\ZM_p (r))$. On obtient ainsi un morphisme de modules galoisiens $U'_\infty(r-1) \longrightarrow H^1(G_S,\ZM_p (r))$ et puisque l'action de ${G_\infty}$ est triviale à droite on obtient un morphisme $\al\colon\U'_\infty(r-1)_{G_\infty} \longrightarrow H^1(G_S,\ZM_p (r))$ caractérisé par les formules (pour tout $m\geq 1$ et tout $\overline{u_\infty\otimes t(r-1)}=\overline{(u_n)_{n\in\NM}\otimes t(r-1)}\in \U'_\infty(r-1)_{G_\infty}$)~:
\begin{equation}\label{forexpal}
(\varphi_m \circ \al) (\overline{u_\infty\otimes t(r-1)}) =
\sum_{g\in G_m} (\kappa^{r-1}(g) g u_{m} )\otimes t(r-1)_m ,
\end{equation}
où l'on identifie $\U'_m \otimes \ZM/p^m (r-1)$ avec son image dans $H^1(G_S(F_m),\ZM/p^m (r))$
par l'homomorphisme de Kummer et où $\varphi_m\colon H^1(G_S,\ZM_p (r)) \longrightarrow H^1(G_S(F_m),\ZM/p^m (r))$ est la composée des applications (commutantes) de réduction modulo $p^m$ et de restriction. On doit aussi remarquer que par définition $G_m=\Gal(F(\ze_{p^m})/F)$ opère sur $\mu_{p^m}$ et en particulier $\kappa(g)$ est bien défini modulo $p^m$ pour tout $g\in G_m$.

\begin{prop} On note $X'_\infty$ la limite de projective pour les applications de normes des $p$-parties des $(p)$-classes $X'_n$ de $F_n$. Le morphisme $\al$ induit la suite exacte
\begin{equation}\label{seknf}
\xymatrix{0\ar[r]& (\U'_\infty(r-1))_{G_\infty} \ar[r]^-\al & H^1(G_S,\ZM_p(r)) \ar[r] & X'_\infty(r-1)^{G_\infty} \ar[r] &0 }.
\end{equation}
\end{prop}
%
\begin{proof}
L'application $\alpha$ coïncide avec celle définie dans la preuve du théorème 3.2 bis (\cite{KNF} p.691).
\end{proof}

\section{\'Eléments cyclotomiques et relations de normes tordues}\label{eltwis}

Le conducteur de $F$ est de la forme $f=p^a d$ avec $p\nmid d$. Pour tout $b\mid d$ et tout $n\geq a$ on note
$$\ep_{b,n}^F=N_{\QM(\ze_{bp^{n}})/F_n \cap \QM(\ze_{bp^{n}})} (1-\ze_{bp^{n}}).$$
Pour $n\geq a$, les $p$-unités $\ep_{b,n}^F$ sont cohérentes en norme suivant $N_{F_{n+1}/F_n}$ et définissent donc des éléments spéciaux $$\ep^F_{b,\infty}=(\ep_{b,n}^F)_{n\geq a} \in \U'_\infty.$$
Et pour obtenir la cohérence en normes pour tout $n$ on convient de noter pour $n<a$~:
$$\ep_{b,n}^F=N_{F_a/F_n}  (\ep_{b,a}^F).$$
Si un nombre premier $\ell\neq p$ divise $d$ et pas $b$ alors il est non ramifié dans $\QM(\ze_{bp^\infty})/\QM$ et le Frobenius $Fr_\ell$ qui agit sur toutes les  racines de l'unité d'ordre premier à $\ell$ par multiplication par $\ell$ est un élément bien défini de l'algèbre de groupe complète
$\ZM_p[[\Gal(\QM(\ze_{bp^\infty})\cap F_\infty /\QM)]]$.
 Les éléments cyclotomiques sont liés par les relations de normes bien connues (pour tout premier $\ell$ et tout entier $b$ avec  $\ell b \mid d$)~:
\begin{equation}\label{relnor}N_{\QM(\ze_{\ell bp^\infty})\cap F_\infty / \QM(\ze_{bp^\infty})\cap F_\infty} (\ep^F_{\ell b,\infty}) = \left \{
\begin{array}{l}
(1-Fr_\ell^{-1}) \ep^F_{b,\infty} \ \mathrm {si}\  \ell\nmid b\\
\ep^F_{b,\infty} \ \mathrm{sinon.} \\
\end{array} \right . \end{equation}
Le personnage principal de notre étude est le codescendu du  $r$-ième tordu de ces éléments cyclotomiques.
\begin{defi} On appelle $r$-ième élément spécial $p$-adique de conducteur $b$ et l'on note $c^F_b(r)\in H^1(G_S,\ZM_p(r))$ l'image par $\al$ de $\overline {\ep^F_{b,\infty}\otimes t(r-1)}\in (\U'_\infty(r-1))_{G_\infty}$~:
\begin{equation}
c^F_b(r):=\al \left (\overline {\ep_{b,\infty}^F\otimes t(r-1)}\right ).
\end{equation}
L'élément spécial attaché à $F$ sur lequel on s'attardera plus particulièrement est celui de conducteur $d$~:
\begin{equation}
c^F(r):=c^F_d(r)=\al \left (\overline {\ep_{d,\infty}^F\otimes t(r-1)}\right ).
\end{equation}
\end{defi}
Ces éléments spéciaux sont caractérisés par leurs images contre les morphismes
$\displaystyle \varphi_m\colon H^1(G_S,\ZM_p(r)) \longrightarrow H^1(G_S(F_m),\ZM/p^{m}(r))$.
\begin{lem} Pour tout $b\mid d$ et tout $m\in \NM$ on a~:
\begin{equation}\label{phi(c)0} \varphi_m(c_b^F(r)) = \sum_{g\in G_m} \kappa^{r-1}(g) g (\ep_{b,m}^F) \otimes t(r-1)_m.
\end{equation}
En particulier pour $m\geq a$ et pour tout $b\mid d$ on a~:
\begin{equation}\label{phi(c)} \varphi_m(c_b^F(r)) = \sum_{g\in \Gal(\QM(\ze_{bp^m})/\QM(\ze_{bp^a})\cap F)} \kappa^{r-1}(g) g (1-\ze_{bp^m}) \otimes t(r-1)_m.
\end{equation}
\end{lem}
\begin{proof} Pour obtenir (\ref{phi(c)0}) il suffit d'appliquer la formule (\ref{forexpal}) à la suite $(u_n)=\ep_{b,\infty}$. Ici aussi les quantités $\kappa(g)$ sont bien définies modulo $p^m$ parce que $\Gal(\QM(\ze_{bp^m})/\QM)$ opère sur $\mu_{p^m}$. Pour en déduire la formule (\ref{phi(c)}) on reprend la définition  $\ep_{b,m}^F=N_{\QM(\ze_{bp^{m}})/F_m \cap \QM(\ze_{bp^{m}})} (1-\ze_{bp^{m}})$. Or $\Gal(\QM(\ze_{bp^{m}})/F_m \cap \QM(\ze_{bp^{m}}))$ agit trivialement sur les racines $p^m$-ièmes de l'unités et donc pour $h\in \Gal(\QM(\ze_{bp^{m}})/F_m \cap \QM(\ze_{bp^{m}}))$ on a $\kappa^{r-1}(h)\equiv 1 [p^m]$. Par identification de $G_m$ avec $\Gal(F_m \cap \QM(\ze_{bp^{m}})/F\cap \QM(\ze_{bp^{m}}))$ et par transitivité de la norme on obtient l'identité (\ref{phi(c)}).\end{proof}

{\it Remarque :} Par construction même, l'image dans $H^1(G_S(\QM(\ze_f)),\ZM_p(r))$ de l'élément spécial $c^F_b(r)$
est l'élément noté $c_r(\ze_f)$ et appelé élément de Soulé-Deligne dans la définition 3.1.2 de \cite{HK03}.
Pour les liens entre ces éléments, ceux définis par Deligne dans \cite{D89}, ceux définis par Soulé dans la section 4.4 de \cite{Sou87} (voir aussi \cite{Sou80}), ceux définis par Beilinson dans \cite{Bei83}, voir les articles originaux et les comparaisons faites dans \cite{HK03} et \cite{HuW}.

\begin{prop} Les $r$-ièmes éléments spéciaux de conducteurs différents sont liés par les relations de normes (pour tout nombre premier $\ell$ et tout $b$ tel que $\ell b\mid d$)~:
\begin{equation}\label{reltor}N_{\QM(\ze_{\ell bp^\infty})\cap F / \QM(\ze_{bp^\infty})\cap F} (c_{\ell b}^F(r)) = \left \{
\begin{array}{l}
(1-\ell^{r-1} Fr_\ell^{-1}) c^F_{b} (r)\ \mathrm {si}\ \ell \nmid b\\
c^F_{b}(r) \ \mathrm{sinon} \\
\end{array} \right . \end{equation}
\end{prop}
\begin{proof} On remarque tout d'abord que $\Gal(\QM(\ze_{\ell bp^\infty})\cap F / \QM(\ze_{bp^\infty})\cap F)$ agit trivialement sur le module de Tate $\ZM_p(r-1)$, tandis que, comme $\ell \neq p$, on a  $Fr_\ell\ t(r-1)=\ell^{r-1} t(r-1)$.
On part de la relation (\ref{relnor}) qui, après twist à la Tate, donne :
\begin{eqnarray}
\nonumber N_{\QM(\ze_{\ell bp^\infty})\cap F / \QM(\ze_{bp^\infty})\cap F} \left (\ep^F_{\ell b,\infty}\otimes t(r-1)\right )  & =
N_{\QM(\ze_{\ell bp^\infty})\cap F / \QM(\ze_{bp^\infty})\cap F} \left (\ep^F_{\ell b,\infty}\right )\otimes t(r-1) \\
\nonumber & =
\left \{
\begin{array}{l}
((1-Fr_\ell^{-1}) \ep^F_{b,\infty}) \otimes t(r-1)\ \mathrm {si} \ \ell\nmid b\\
\ep^F_{b,\infty}\otimes t(r-1) \ \mathrm{sinon.} \\
\end{array} \right .
\end{eqnarray}
Le cas particulier $\ell\mid b$ de la formule (\ref{reltor}) s'en déduit alors immédiatement en prenant les co-invariants puis en appliquant $\al$. Lorsque $\ell\nmid b$ on a :
$$\begin{array}{rl}
 N_{\QM(\ze_{\ell bp^\infty})\cap F / \QM(\ze_{bp^\infty})\cap F} & \left (\ep^F_{\ell b,\infty}\otimes t(r-1)\right )  =
((1-Fr_\ell^{-1}) \ep^F_{b,\infty}) \otimes t(r-1) \\
 & =  \left (\ep^F_{b,\infty} \otimes t(r-1) \right ) -
\left ( Fr_\ell^{-1} \ep^F_{b,\infty}  \otimes \ell^{r-1} Fr_\ell^{-1}  t(r-1)\right )
\\ &= \left ( 1-\ell^{r-1} Fr_\ell^{-1}\right ) \ep^F_{b,\infty} \otimes t(r-1).
\end{array}$$
Et on conclut à nouveau en prenant les co-invariants et en appliquant $\al$.
\end{proof}
\begin{defi}\label{CFr} On appelle groupe des éléments cyclotomiques et on note $C^F(r)$ le sous-module de $H^1(G_S,\ZM_p(r))$ engendré par les $r$-ièmes éléments spéciaux de conducteurs divisant $d$~:
$$C^F(r)=\langle c_b^F(r),\, b\mid d \rangle $$
\end{defi}

L'un des intérêts de ce groupe réside dans l'analogue supérieur de la formule analytique du nombre de classes démontrée par Kolster, Nguy{\cftil{e}}n-Quang-{D}{\cftil{o}} et Fleckinger. Par exemple on a l'égalité~:
\begin{theo}\label{knfmain} On rappelle que $F$ est totalement réel et que $r> 1$ est impair.
$$\# \left ( \frac {H^1(G_S,\ZM_p(r))} {C^F(r)} \right ) =\# H^2(G_S,\ZM_p(r))$$
\end{theo}
\begin{proof} Soit $G$ un groupe abélien agissant sur un corps de nombre $K$. Pour tout caractère $\chi$ sur $G$, si $p\nmid |G|$ on peut définir les $\chi$-composantes avec des idempotents de $\ZM_p[G]$.  Toujours avec l'hypothèse $p\nmid |G|$, le passage à ces $\chi$-composantes est exact. D'après \cite{KNF} theorem 5.4, si $K$ est une extension imaginaire, absolument abélienne de $F$ tel que $p \nmid [K:F]$ on a pour tout $\chi$ sur $G$ tel que $\chi(-1)=(-1)^{r-1}$~:  
\[\# \left ( \frac {H^1(G_S(K),\ZM_p(r))} {C^K(r)} \right )^\chi =\# H^2(G_S(K),\ZM_p(r))^\chi.\]
(Noter la différence de notation entre \cite{KNF} et le présent article, les rôles de $K$ et $F$ sont inversés et la notation $\Cba'(r-1)_{G_\infty}$ est ici allégée en $C^K(r)$.) Maintenant pour démontrer la formule du théorème \ref{knfmain} il suffit de prendre $K=F[i]$ et comme $p\neq 2$ le théorème s'applique avec $\chi$ le caractère trivial sur $\Gal(K/F)$ et $r>1$ impair.\end{proof}
La philosophie générale de la $K$-théorie arithmétique prédit que les $H^1(G_S,\ZM_p(r))$ pour les $r$ impairs se comportent comme des $p$-unités tordues, cette intuition étant précisée notament par la suite exacte (\ref{seknf}); tandis que les $H^2(G_S,\ZM_p(r))$ se comportent comme des $p$-groupes de $S$-classes d'idéaux tordues, voir notamment le lemme 4.3 p. 269 de \cite{Kur99}. De ce point de vue l'égalité du théorème \ref{knfmain} est l'analogue supérieur et sans facteur parasite des formules d'indices de Sinnott pour le twist $r=1$ (\cite{Si78,Si80}).

{\it Remarque :} Lorsque $p\nmid [F:\QM]$ la démarche de \cite{KNF} s'applique aussi pour démontrer la formule 
$$\# \left ( \frac {H^1(G_S,\ZM_p(r))} {C^F(r)} \right )^\chi =\# H^2(G_S,\ZM_p(r))^\chi,$$
pour tout caractère $\chi$ sur $F$. Le referee demande que soit précisé l'état des connaissances sur l'égalité entre les ordres des $\chi$-parties (respectivement des $\chi$-quotients) des modules ci-dessus lorsque $p$ divise le degré $[F:\QM]$. C'est une question difficile, dont la réponse dépend a priori du choix entre $\chi$-parties 
et $\chi$-quotient, et n'est pas connue des auteurs.  

Il est de toute façon intéressant de trouver des méthodes de calcul explicite des ordres des $\chi$-composantes du quotient
$H^1(G_S,\ZM_p(r))/C^F(r)$ pour \underline{tous} les caractères $\chi$ de $F$. C'est ce qui a été fait (parmi d'autres résultats) dans l'article \cite{Kur99} mais en supposant que l'exposant de $G$ divise $p-1$, et en particulier que $p\nmid [F:\QM]$. Dans le cas général, traité ici, il faut aussi choisir une notion précise de $\chi$-composantes. En effet au moins deux foncteurs naturels donnant des $\chi$-composantes co-existent. On reviendra sur ces $\chi$-composantes en section \ref{chipart}. Dans la section qui suit, pour mieux décrire les composantes auxquelles on s'intéresse, on calcule le caractère de quelques représentations linéaires.

\section{Caractères des représentations en jeu}\label{carut}
\begin{defi} Soit $G$ un groupe fini.
\begin{enumerate}
\item  Pour tout $\ZM_p[G]$ module $M$, on appelle caractère de $M$ la fonction centrale sur $G$
obtenue en prenant la trace de $M\otimes \QM_p$ vue comme représentation de $G$.
\item On appelle caractère régulier sur $G$ le caractère de $\ZM_p[G]$. On appelle caractère trivial sur $G$ le caractère de $\ZM_p$ avec action triviale de $G$.
\end{enumerate}
\end{defi}
Il est bien connu que deux $\ZM_p[G]$-modules deviennent isomorphes après extension des scalaires de $\ZM_p$ à $\QM_p$ si et seulement si ils ont même caractères.
\begin{prop}\label{car} On rappelle que le twist $r\neq 1$ est supposé impair et le corps de base $F$ totalement réel.
\begin{enumerate}
\item $C^F(r)$ est d'indice fini dans $H^1(G_S,\ZM_p(r))$.
\item $\al((\U'_\infty(r-1))_{G_\infty})$ est d'indice fini dans $H^1(G_S,\ZM_p(r))$.
\item Le caractère commun sur $\Gal(F/\QM)$ des trois modules $C^F(r),\al((\U'_\infty(r-1))_{G_\infty})$ et $H^1(G_S,\ZM_p(r))$ est le caractère régulier.
\end{enumerate}
\end{prop}

\begin{proof} 1-- est une conséquence du théorème 5.4 de \cite{KNF} et 2-- de la suite exacte (\ref{seknf}).
Pour 3-- il suffit de calculer le caractère de $(\U'_\infty(r-1))_{G_\infty}$ par injectivité de $\al$. Parce que $r$ est impair et $F$ totalement réel, ce caractère est le caractère régulier, autrement dit $(\U'_\infty(r-1))_{G_\infty}\otimes \QM_p \cong \QM_p[G]$. Ce résultat fait probablement partie du folklore et par exemple est compatible avec les rangs mentionnés dans \cite{Kosur}. Nous n'avons malheureusement pas trouvé de référence qui précise le caractère des modules qui interviennent dans cet article. Pour la commodité du lecteur on va énoncer et refaire la preuve du lemme à suivre qui conclut la démonstration de la proposition \ref{car}.
\end{proof}

\begin{lem}\label{caruinf} Soit $\UC_\infty$ (resp. $\U_\infty$) la limite projective des unités semi-locales (resp. globales) le long de la tour $F(\ze_{p^\infty})/F$, relativement aux applications de norme.
\begin{enumerate}
\item Pour tout entier $r\neq -1$ le caractère de $\UC_\infty(r)_{G_\infty}$ est le caractère régulier. Le caractère de $\UC_\infty(-1)_{G_\infty}$ est la somme du caractère régulier et de celui de $\ZM_p[S]$.
\item Si $r$ est pair alors le module $\U_\infty(r)_{G_\infty}$ est de caractère régulier. Le caractère de $\U_\infty(-1)_{G_\infty}$ est le caractère de $\ZM_p$ avec action triviale de $G$. Si $r$ est impair et différent de $-1$ alors $\U_\infty(r)_{G_\infty}$ est de $\ZM_p$-torsion.
\item Les modules $\U_\infty(r)_{G_\infty}$ et $\U'_\infty(r)_{G_\infty}$ ont même caractère.
\end{enumerate}
\end{lem}
\begin{proof}
On note $\GC=\Gal(F_\infty/\QM)$ et $\Lambda[\GC]$ l'algèbre de groupe pro-$p$-complété de $\GC$. Clairement on a $\Lambda[\GC]\cong \Lambda[\GC](r)$ pour tout $r$. Lorsque $p$ est modérément ramifié dans $F$, la théorie de Coleman donne une suite exacte (voir \cite{Grei92}
ou \cite{T99} pour les détails)
\begin{equation}\label{secol}
\xymatrix {0\ar[r] & \oplus_{v\mid p} \ZM_p(r+1) \ar[r] & \UC_\infty(r) \ar[r] & \La[\GC] \ar[r] &  \oplus_{v\mid p} \ZM_p(r+1) \ar[r] & 0 \ .}
\end{equation}
Sans aucune condition de ramification mais si $F$ est abélien sur $\QM$ il existe un sous-corps $L\subset F_\infty$ qui est abélien  sur $\QM$, dans lequel $p$ est modérément ramifié, et tel que $F(\ze_{p^\infty})=L(\ze_{p^\infty})$ (voir par exemple \cite{CJM} lemme 1.2). Ainsi la suite exacte (\ref{secol}) est valide pour tout corps $F$ abélien.

On vérifie le cas $r=-1$ du 1- du lemme \ref{caruinf}. Pour cela on note $\VC_\infty(-1) \subset \La[\GC]$ l'image de $\UC_\infty(-1)$. La suite (\ref{secol}) implique alors l'égalité $\VC_\infty(-1)^{G_\infty}=0$ et redonne par descente les deux suites exactes~:
$$\xymatrix{0 \ar[r] & \ZM_p[S]\ar[r] & \UC_\infty(-1)_{G_\infty} \ar[r] & \VC_\infty(-1)_{G_\infty} \ar[r] & 0 .}$$
$$\xymatrix{0 \ar[r] & \ZM_p[S] \ar[r] & \VC_\infty(-1)_{G_\infty} \ar[r] & \ZM_p[G] \ar[r] & \ZM_p[S] \ar[r] & 0 .}$$
Comme les caractères sont additifs en suites exactes cela démontre la cas particulier $r=-1$ du point 1.
Pour $r\neq -1$ les $G_\infty$ invariants et co-invariants du module $\oplus_{v\mid p} \ZM_p(r+1)$ sont finis et donc leur caractère est nul. La même chasse au diagramme que ci-dessus mais en plus simple démontre donc que $\UC_\infty(r)_{G_\infty}$ et $\La[\GC]_{G_\infty}$ ont même caractère, c'est-à-dire le caractère régulier. Le point 1 est démontré.

Pour 2 on sait, comme $p\neq 2$, que $\U_\infty \cong \ZM_p(1)\oplus \U_\infty^+$. Si $r$ est pair alors la conjugaison complexe de $G_\infty$ agit trivialement sur $\ZM_p(r)$ donc elle agit trivialement sur $\UC_\infty^+(r)$ et agit par $-1$ sur $\UC_\infty^-(r)$. On en déduit $\UC_\infty(r)_ {G_\infty}=(\UC_\infty^+\oplus \UC_\infty^-)(r)_{G_\infty}= \UC_\infty^+(r)_{G_\infty}$. Par la conjecture faible de Leopoldt, qui est un théorème pour l'extension cyclotomique, on sait que $\UC_\infty^+/\U_\infty^+$ est un module de torsion de type fini sur l'algèbre d'Iwasawa classique $\Lambda=\Lambda[\Gal(F_\infty/F(\ze_p))]$. Donc $(\UC_\infty^+/\U_\infty^+)(r)$ aussi et $(\UC_\infty^+/\U_\infty^+)(r)\otimes \QM_p$ est de dimension finie sur $\QM_p$. En conséquence les modules $(\UC_\infty^+/\U_\infty^+)(r)^{G_\infty}$ et $(\UC_\infty^+/\U_\infty^+)(r)_{G_\infty}$ ont même caractère. Dans cette parité de $r$ le point 2 se déduit donc du point 1 et de la suite exacte
$$\xymatrix{0\ar[r] & (\UC_\infty^+/\U_\infty^+)(r)^{G_\infty} \ar[r] & \U_\infty^+(r)_{G_\infty} \ar[r] & \UC_\infty^+(r)_{G_\infty} \ar[r] & (\UC_\infty^+/\U_\infty^+)_{G_\infty}\ar[r] & 0}.$$

On passe au cas où $r$ est impair. Comme $F$ est totalement réel la conjugaison complexe est un élément bien défini de $G_\infty$ qui agit simultanément par $-1$ et par $1$ sur $\U_\infty^+(r)_{G_\infty}$. Ce dernier module est donc trivial. Cela démontre le cas $r$ impair du point 2, puisque $\ZM_p(r)_{G_\infty}$ a bien le caractère requis suivant $r=-1$ ou $r\neq -1$.

On démontre le point 3.
Les modules $\U_\infty$ et $\U'_\infty$ ont la même $\Lambda$-torsion (à savoir $\ZM_p(1)$) et donc $\U'_\infty(r)$ et $\U_\infty(r)$ ont les même $G_\infty$ invariants (à savoir $0$ si $r\neq -1$ et $\ZM_p(0)$ si $r=-1$.)
En utilisant alors que $\U'_\infty(r)/\U_\infty(r)$ s'injecte dans $\ZM_p[S](r)$ qui est de $\La$-torsion on démontre le point 3 en suivant le même type de chasse au diagramme que précédemment.

\end{proof}

On  énonce un premier résultat utile à notre calcul de $\chi$-indice~:

\begin{theo}\label{freespe} L'élément spécial $c^F(r)$ engendre un sous-$\ZM_p[G]$-module libre de rang $1$ de $H^1(G_S,\ZM_p(r))$.
\end{theo}
\begin{proof} Comme on étudie le module monogène engendré par $c^F(r)$, il suffit de montrer que pour tout caractère de Dirichlet $\chi$ sur $G$ l'élément $c^F(r)\otimes e_\chi \in C^F(r)\otimes \CM_p$ est non nul, où $e_\chi$ désigne l'idempotent usuel $e_\chi =(1/\#G) \sum_{g} \chi(g) g^{-1}$. On fixe $\chi$ un tel caractère, on note $F_\chi$ le sous-corps de $F$ fixé par $\ker \chi$ et on note $kp^s$ le conducteur commun à $\chi$ et à $F_\chi$ (avec $s\leq a$). Le sous-espace $\chi$-isotypique  $C^F(r)\otimes e_\chi \subset C^F(r)\otimes \CM_p$ est engendré par les
$ c^F_b(r) \otimes e_\chi$. Si $\QM(\ze_{bp^\infty})\cap F$ ne contient pas
$\QM(\ze_{kp^s})\cap F$ alors $\QM(\ze_{bp^\infty})\cap F$ ne contient pas non plus $F_\chi$ parce que $kp^s$ est le conducteur de $\chi$. Et dans ce cas $c^F_b(r) \otimes e_\chi$ est nul parce que tout élément non trivial $g$ de $\Gal(F_\chi/\QM(\ze_{bp^\infty})\cap F_\chi)$ agit simultanément par multiplication par $1$ et par $\chi(g)\neq 1$ sur $c^F_b(r) \otimes e_\chi$.
Par contre tous les $c^F_b(r) \otimes e_\chi$ pour tous les $b$ tels que $\QM(\ze_{bp^\infty})$ contienne $F_\chi$ sont co-linéaires à $c^F_k(r) \otimes e_\chi$. En effet la formule  (\ref{reltor}) donne~:
$$\begin{array}{rl} [\QM(\ze_{bp^\infty})\cap F : \QM(\ze_{kp^\infty}) \cap F ] c^F_b(r) \otimes e_\chi & = N_{\QM(\ze_{bp^\infty})\cap F / \QM(\ze_{kp^\infty}) \cap F} (c^F_b(r) \otimes e_\chi)\\ &=\prod_{\ell \mid b,\ell \nmid k} (1-\ell^{r-1} Fr_\ell^{-1}) c^F_k(r)\otimes e_\chi, \end{array}$$
c'est-à-dire~:
$$c^F_b(r) \otimes e_\chi = \frac {\prod_{\ell \mid b,\ell \nmid k} (1-\ell^{r-1} \chi(\ell)^{-1})} {[\QM(\ze_{bp^\infty})\cap F : \QM(\ze_{kp^\infty}) \cap F ]}  c^F_k(r)\otimes e_\chi .$$
Pour $r\neq 1$ les facteurs eulériens $(1-\ell^{r-1} \chi(\ell)^{-1})$ ne s'annulent jamais et par la proposition \ref{car} il suit $c^F_f(r)\otimes e_\chi \neq 0$.
\end{proof}

\section{Composantes isotypiques}\label{chipart}
Dans cette section on rappelle les définitions des $\chi$-parties et $\chi$-quotients et l'on décrit précisément quel $\chi$-indice est calculé dans la suite en expliquant nos choix.
\begin{defi} Soit $M$ un $\ZM_p[G]$-module et $\chi$ un caractère sur $G$. Suivant Tsuji \cite{T99}, on note $\underline{\ZM_p[\chi]}$ le $\ZM_p[\chi]$-module libre de rang $1$ muni de l'action de $G$ par multiplication par $\chi$.
\begin{enumerate}
\item On appelle $\chi$-partie de $M$ et on note $M^\chi$ le plus grand sous-module de $M$ sur lequel l'action de $G$ est décrite par $\chi$. Formellement on peut définir $M^\chi$ par
    $$\Hom_{\ZM_p[G]} (\underline{\ZM_p[\chi]},M)\cong M^\chi \subset M.$$
\item On appelle $\chi$-quotient de $M$ et on note $M_\chi$ le plus grand quotient de $M$ sur lequel l'action de $G$ est décrite par $\chi$. Formellement on peut définir $M_\chi$ par
    $$M\twoheadrightarrow M_\chi\cong M\otimes_{\ZM_p[G]} \underline{\ZM_p[\chi]}.$$
\end{enumerate}
\end{defi}
Pour les propriétés basiques de ces deux foncteurs on pourra consulter le \S II de \cite{So90}, la section 2 de \cite{T99}, ou la thèse du premier auteur \cite{Tan}. Il serait beau d'avoir des méthodes simples pour calculer les ordres $\#(H^1(G_S,\ZM_p(r))/C^F(r))^\chi$ ou $\#(H^1(G_S,\ZM_p(r))/C^F(r))_\chi$ qui en général ne coïncident pas forcément.
Dans la suite on choisit de calculer $\#(H^1(G_S,\ZM_p(r))^\chi/\langle c^F(r) \rangle^\chi)$, qui se prête le mieux à l'approche de cet article, par exemple en raison du théorème \ref{freespe}. Très souvent (en particulier par exemple si $p\nmid [F:\QM]$) tous ces $\chi$-indices sont égaux. Mais dans la généralité avec laquelle cet article est rédigé il semble utile d'expliquer ce choix. Dans le but d'étudier l'arithmétique des valeurs spéciales d'une fonction $L$ de Dirichlet (primitive) $L(s,\chi)$ il est naturel de fixer un caractère de Dirichlet $\chi$ de conducteur $f$, de prendre $F=F_\chi:= \QM(\ze_f)^{\ker \chi}$ et d'utiliser le module monogène engendré par l'élément cyclotomique $\langle c^F(r)\rangle$ à la place du module $C^F(r)$. Dans cette restriction (naturelle) en généralité, on a alors $G$ cyclique, et notre choix n'est pas restrictif à cause de la proposition suivante~:
\begin{prop} On suppose $G$ cyclique.
\begin{enumerate}
\item Pour tout $\ZM_p[G]$-module fini $M$ on a $\# M^\chi=\#M_\chi$, en particulier on a l'égalité
$$ \#\displaystyle \left (\frac {H^1(G_S,\ZM_p(r))}{\langle c^F(r)\rangle}\right )^\chi=\# \displaystyle \left (\frac {H^1(G_S,\ZM_p(r))}{\langle c^F(r)\rangle}\right )_\chi\ .$$
\item On a un isomorphisme $$\displaystyle \frac{H^1(G_S,\ZM_p(r))^\chi} {\langle c^F(r) \rangle^\chi}\overset \sim \longrightarrow \left (\frac {H^1(G_S,\ZM_p(r))}{\langle c^F(r)\rangle}\right )^\chi.$$
\end{enumerate}
\end{prop}
\begin{proof}
Lorsque $G$ est cyclique engendré par $g$ on a une suite exacte
$$\xymatrix{ 0 \ar[r] & M^\chi \ar[r] & M\otimes \ZM_p[\chi]  \ar[r]^-{g-\chi(g)} &   M\otimes \ZM_p[\chi]
\ar[r] & M_\chi \ar[r] & 0 }$$
Cela donne l'égalité des ordres du point 1.

Par le théorème \ref{freespe} le module $\langle c^F(r)\rangle$ est isomorphe à $\ZM_p[G]$ et donc on a les isomorphismes
$\langle c^F(r)\rangle_\chi \cong \ZM_p[\chi] \cong \langle c^F(r)\rangle^\chi$. Par le lemme 2.2 de \cite{T99} on a une suite exacte longue~:
$$\xymatrix@=18pt{0\ar[r] &  \langle c^F(r)\rangle^\chi \ar[r] & H^1(G_S,\ZM_p(r))^\chi \ar[r] & \left (\frac {\displaystyle H^1(G_S,\ZM_p(r))}{\displaystyle \langle c^F(r)\rangle} \right)^\chi \ar[r]^-\de & \langle c^F(r)\rangle_\chi \ar[r] & \cdots }.$$
Mais le module $\langle c^F(r)\rangle_\chi$ est sans $\ZM_p$-torsion tandis que $(H^1(G_S,\ZM_p(r))/\langle c^F(r)\rangle)^\chi$ est fini. Donc la flèche de connexion $\de$ est triviale et cela donne l'isomorphisme du point 2.
\end{proof}
Comme cela ne présente pas de difficultés supplémentaires, dans la suite on prend $F$ abélien totalement réel quelconque, on prend $\chi$ un caractère sur $G$ quelconque et l'on va calculer l'indice $\#(H^1(G_S,\ZM_p(r))^\chi/\langle c^F(r) \rangle^\chi)$ en utilisant les morphismes de réduction modulo $\ell$ dans le style de Kurihara \cite{Kur99}.

\section{Les morphismes de réduction modulo $\ell$}\label{redl}
Soit un entier $n\geq 1$ et étant donné un nombre premier $\ell\in \ZM$, tel que $\ell\equiv 1[dp^n]$, et un idéal fixé $\la$ de $F$ divisant $\ell$, l'objet de cette section est de construire explicitement et de décrire des propriétés d'une application de réduction $G$-équivariante~:
\begin{equation}\label{phi}
\phi_{\ell,n} \colon H^1(G_S,\ZM/p^n(r)) \longrightarrow \bigoplus_{g\in G} H^1(\FM_{\la^g},\ZM/p^n(r))
\end{equation}
où $\FM_\la$ désigne le corps résiduel $O_F/\la\cong \FM_\ell$. Cette application $\phi_{\ell,n}$ est essentiellement "l'application naturelle" de \cite{Kur99} p. 265, mais pour l'usage ultérieur fait dans cet article on l'examine un peu plus en détail ici. Dans l'article \cite{Ga01} Gajda examine aussi une application naturelle en $K$-théorie
$\phi_{p,l}\colon K_{2n+1}(\ZM) \longrightarrow K_{2n+1} (\FM_p)[l^\infty]$. Sous la conjecture de Quillen-Lichtenbaum, via les classes de Chern, le $\phi_{l,p}$ de Gajda  donne ainsi une application naturelle de réduction pour le cas particulier où le corps $F$ est un sous-corps de $\QM(\ze_p)$. En supposant l'hypothèse de Vandiver, Gajda démontre que son $\phi_{l,p}$ est non nul pour un ensemble de premiers $p$ dont il calcule la densité. Notre démarche requiert des informations plus précises que la non trivialité des $\phi_{\ell,n}$ et se situe en cohomologie Galoisienne.

On rappelle que $\GC_S$ désigne $\Gal(\Omega_S/\QM)$.
Pour tout idéal premier $\la$ de $O_F$ divisant $\ell\equiv 1[dp]$, on se donne un idéal $\LC=\LC(\la)$ de $\Omega_S$ divisant $\la$. Soit
$D_\LC\subset G_S$ le sous-groupe de décomposition dans $\GC_S$ de $\LC$. Alors pour tout $p$-groupe discret $M$ muni d'une action continue de $\GC_S$ on a un morphisme de restriction
$res_\la \colon H^1(G_S,M) \longrightarrow H^1(D_\LC,M)$. Si on se donne en outre $\si\in \GC_S$, alors cette restriction commute avec le morphisme de conjugaison $\si_*$ qui définit l'action de $\si$ sur les $H^1(H,M)$ pour tout sous-groupe fermé $H$ de $\GC_S$ (voir par exemple \cite{NSW} p. 44). Ainsi modulo l'identification $D_\LC \overset \si \cong D_{\LC^\si}$ ce morphisme de restriction dépend seulement de $\la$ et pas du choix de $\LC$ divisant $\la$, parce qu'un autre choix conjugue $\LC$ par un $\si \in G_S$ qui agit trivialement sur $H^1(G_S,M)$. Par contre pour un $\si\in \GC_S$ quelconque qui définit un idéal $\la^\si$ de $F$ divisant $\ell$ on a un diagramme commutatif
\begin{equation}\label{Geq}
\xymatrix {H^1(G_S,M) \ar[r]^-{res_\la}  \ar@{^{(}->>}[d]^-{\si_*} & H^1(D_\LC, M ) \ar@{^{(}->>}[d]^-{\si_*} \\
    H^1(G_S,M) \ar[r]^-{res_{\la^\si}}  & H^1(D_{\LC^\si}, M ) .}
\end{equation}
Si $\la\not\in S$, alors $\la$ est non ramifié dans $\Omega_S$ et le pro-$p$-quotient maximal de $D_\LC$ s'identifie avec le pro-$p$-quotient maximal de $\Gal(\overline{\FM}_\ell/\FM_\la)$. Donc pour tout $p$-groupe discret $M$ on a des  isomorphismes $H^1(D_{\LC^\si},M)  \cong H^1(\FM_{\la^\si},M)$. Ceci posé on construit un morphisme $\phi_{\ell,n}$ comme annoncé dans le diagramme (\ref{phi}) en prenant $M=\ZM/p^n(r)$ dans ce qui précède et en posant $$\phi_{\ell,n}=\bigoplus_{g\in G}  res_{\la^g}.$$  On obtient alors immédiatement que $\phi_{\ell,n}$ est $G$-équivariant pour l'action par permutation de $G$ sur la somme directe $\displaystyle \bigoplus_{g\in G} H^1(\FM_{\la^g},\ZM/p^n(r))$.

Pour rendre plus explicite cette application $\phi_{\ell,n}$ on va en donner une interprétation kummerienne.
\begin{defi} Soit $\la\subset O_F$ divisant $\ell\equiv 1[dp^n]$, $\ell\neq 2$ un nombre premier et $\LC$ le premier fixé de $\Omega_S$ divisant $\la$. Par abus on écrit encore $\LC$ pour $\LC\cap F_n$, et on identifie $O_{F_n}/\LC$ avec $\FM_\la$.
\begin{enumerate}
\item On note $red_\la\colon F_n^\times/p^n \longrightarrow \FM_\la^\times/p^n$ l'unique morphisme de groupe qui fasse commuter le diagramme
\begin{equation*}
\xymatrix{F_n^\times \ar[r] \ar[d] & F_{n,\LC}^\times \ar[r]^-\sim & \FM_\la^\times \times (1+\ell)^{\ZM_\ell} \times \ell^\ZM \ar[r] &
\FM_\la^\times\ar[d]  \\
F_n^\times/p^n \ar[rrr]^-{red_\la} & & & \FM_\la^\times/p^n \\}
\end{equation*}
Explicitement $red_\la$ associe à la classe modulo les puissances $p^n$-ièmes de $x\in F_n^\times$ la classe modulo $\LC$ et les puissances $p^n$-ièmes de $x\ell^{-v_\LC(x)}$.
\item Pour $r\geq 0$ on note $red_\la(r) \colon F_n^\times/p^n (r) \longrightarrow  \FM_\la^\times/p^n (r)$ le $r^{\text{ième}}$ tordu de $red_\la$. Explicitement $red_\la(r)$ envoie $\xba \otimes t(r)_n$ sur $(x\ell^{-v_\LC(x)} + \LC) \otimes t(r)_n$.
\end{enumerate}
\end{defi}
Cette définition permet de calculer explicitement les images de $res_\la$ en utilisant le diagramme commutatif suivant~:
\begin{equation}\label{kumres}
\xymatrix{(F_n^\times/p^n) (r-1) \ar[rd]^-{red_\la(r-1)} \\
 \{x \in F_n^\times/p^n;\; \forall v\nmid p, \, v(x)\equiv 0 [p^n]\} (r-1)\ar@{^{(}->>}[d] \ar[r] \ar@{_{(}->}[u]& (\FM_\la^\times/p^n) (r-1) \ar@{^{(}->>}[d]\\
H^1(G_S(F_n),\ZM/p^n(1))(r-1) \ar[r] & H^1(D_\LC,\ZM/p^n (1)) (r-1) \\
H^1(G_S,\ZM/p^n(r))\ar[u] \ar[ru]_-{res_\la} & }
\end{equation}
Les applications du triangle du bas sont des restrictions; les deux triangles commutent de façon tautologique. Les isomorphismes verticaux s'obtiennent en utilisant le  cocycle de Kummer. La commutativité du carré central est un simple jeu d'écriture à partir de ce cocycle (et il suffit de vérifier cette commutativité pour $r=1$).

Le résultat principal de cette section est le théorème d'isomorphisme suivant qui permet d'utiliser les morphismes de réduction pour calculer les $\chi$-indices~:
\begin{theo}\label{inj} Soit $\chi$ un caractère de $\Gal(F/\QM)$. On Rappelle que $f$, le conducteur de $F$, vaut $f=dp^a$ avec $p\nmid d$. On prend un $n>a$. Il existe un (une infinité de) nombre(s) premier(s) $\ell\equiv 1 [dp^n]$ tel(s) que le morphisme $\phi_{\ell,n}$ induise un isomorphisme entre les $\chi$-parties~:
$$H^1(G_S,\ZM_p(r))^\chi/p^n H^1(G_S,\ZM_p(r))^\chi \overset \sim \longrightarrow \left ( \bigoplus_{g\in G}  \FM_{\la^g}^\times/p^n (r-1)\right )^\chi.$$
\end{theo}
\begin{proof} Tout d'abord les deux modules sont finis, isomorphes sur $\ZM_p[G]$ à $\ZM_p[\chi]/p^n$. Pour $H^1(G_S,\ZM_p(r))$ c'est parce qu'il  est sans $\ZM_p$-torsion (lemme \ref{torfree}) et de caractère régulier (proposition \ref{car}); et donc la $\chi$-partie $H^1(G_S,\ZM_p(r))^\chi$ est un $\ZM_p[\chi]$-module sans torsion de rang $1$. Pour le module de droite $ \left ( \bigoplus_{g\in G}  \FM_{\la^g}^\times/p^n (r-1)\right )^\chi\cong (\ZM/p^n [G])^\chi$, c'est évident.
Il suffit donc de trouver  un (une infinité de)  $\ell$ tels que l'application $$H^1(G_S,\ZM_p(r))^\chi/p^n H^1(G_S,\ZM_p(r))^\chi\longrightarrow \left ( \bigoplus_{g\in G}  \FM_{\la^g}^\times/p^n(r-1)\right )^\chi$$ déduite  de $\phi_{\ell,n}$ soit injective.
On a une suite d'injection~:
$$ H^1(G_S,\ZM_p(r))/p^n \overset {(1)} \hookrightarrow H^1(G_S,\ZM/p^n (r)) \overset{(2)} \hookrightarrow
H^1(G_S(F_n),\ZM/p^n (r))\overset {(3)} \hookrightarrow F_n^\times/p^n (r-1).$$
L'injectivité du morphisme $(1)$ suit de la cohomologie de la suite
$$0 \longrightarrow \ZM_p (r) \overset {p^n}  \longrightarrow \ZM_p(r) \longrightarrow \ZM/p^n (r) \longrightarrow 0.$$ L'injectivité du morphisme $(2)$ suit par la suite d'inflation restriction grâce à la $\Gal(F_n/F)$-trivialité cohomologique des $\ZM/p^n (r)$ pour $r\neq 0$ (lemme de Tate). L'injectivité du morphisme $(3)$ vient de la théorie de Kummer et a déjà été affirmée dans le diagramme (\ref{kumres}). Soit $\LM_n$ l'ensemble des premiers rationnels $\ell\equiv 1 [dp^n]$ et tels que $p^{n+1}\nmid (\ell-1)$. Par le théorème de Dirichlet $\LM_n$ est infini.
On considère l'application de réduction $red(r-1)$ produit des $red_\la(r-1)$ pour $\la$ parcourant les premiers de $F$ divisant les $\ell\in \LM_n$~:
$$red(r-1) \colon F_n^\times/p^n (r-1)\longrightarrow  \prod_{\ell \in \LM_n} \left (\prod_{\la \mid \ell, \la\subset F} \FM_\la^\times/p^n(r-1) \right )$$
\begin{lem} Le morphisme $red(r-1)$ est injectif.
\end{lem}
\begin{proof} Le foncteur $\otimes \ZM/p^n(r-1)$ est exact sur les $\ZM/p^n$-modules, il suffit donc de vérifier cette injectivité pour $r=1$. Soit $\varphi\colon F_n^\times\longrightarrow \prod_{\ell \in \LM_n} \left (\prod_{\la \mid \ell, \la\subset F} \FM_\la^\times \right )$ l'application de réduction avant de prendre le quotient modulo $p^n$. Si $x$ est dans le noyau de $\varphi$ alors $(x-1)$ appartient à tous les $\LC$ avec $v_{\LC}(x)=0$ soit une infinité d'idéaux premiers distincts, et donc $x=1$. Cela montre l'injectivité de $\varphi$. On en déduit un morphisme injectif $$\imath\colon F_n^\times/p^n \longrightarrow \left ( \prod_{\ell \in \LM_n} \prod_{\la \mid \ell, \la\subset F} \FM_\la^\times \right )/ \varphi (p^n F_n^\times).$$ L'image de $\imath$ est contenu dans les éléments d'ordre une puissance de $p$ du quotient de droite. Ce sous-groupe des éléments d'ordre une puissance de $p$ s'injecte lui-même naturellement dans
$ \prod_{\ell \in \LM_n}\left (\prod_{\la \mid \ell, \la\subset F} \FM_\la^\times/p^n \right )$ parce que $p^n$ divise exactement l'ordre $\ell-1$ des $\FM_\la^\times$, pour $\ell\in \LM_n$.
\end{proof}
On reprend la preuve du théorème \ref{inj}. Par exactitude à gauche du foncteur $\chi$-parties on déduit de ce qui précède une suite de morphismes injectifs~:
$$H^1(G_S,\ZM_p(r))^\chi/p^n  \hookrightarrow  \left (H^1(G_S,\ZM_p(r))/p^n\right )^\chi \hookrightarrow \prod_{\ell \in \LM_n}\left (\prod_{\la \mid \ell, \la\subset F_n} \FM_\la^\times/p^n(r-1) \right )^\chi .$$
Pour un $\ell$ fixé, l'application déduite par $\chi$-partie de $\phi_{\ell,n}$ est la composée de cette suite d'injection avec la projection sur $\left (\prod_{\la \mid \ell, \la\subset F_n} \FM_\la^\times/p^n(r-1) \right )^\chi$~: c'est une conséquence de la commutativité du diagramme (\ref{kumres}) (compatibilité entre $res_\la$ et $red_\la(r-1)$).
On fixe $x$ un $\ZM_p[\chi]$-générateur du $\ZM_p[\chi]$-module libre $H^1(G_S,\ZM_p(r))^\chi$. Alors l'image de $x$
dans $H^1(G_S,\ZM_p(r))^\chi/p^n$ et donc aussi
dans $\prod_{\ell \in \LM_n}\left (\prod_{\la \mid \ell, \la\subset F_n} \FM_\la^\times/p^n(r-1) \right )^\chi$ engendre un module libre de rang $1$ sur $\ZM_p[\chi]/p^n$. Il en est forcément de même pour au moins une des projections  $\phi_{\ell,n}^\chi(x)$ et donc pour au moins un $\ell \in \LM_n$. En outre le même raisonnement s'applique si on remplace $\LM_n$ par n'importe lequel de ses sous-ensemble infini $\EC\subset \LM_n$. Donc il existe aussi une infinité de $\ell\in \LM_n$ tel que $\phi_{\ell,n}^\chi$ soit un isomorphisme : cela démontre le théorème \ref{inj}.
\end{proof}

\section{Le théorème d'indice}\label{thmind}

Soit $n\in \NM$ et $\ell$ un nombre premier tel que $dp^n$ divise $\ell-1$. Soit $\chi$ un caractère sur $G$. Dans la section \ref{redl} on a étudié le morphisme
$$\phi_{\ell,n}^\chi \colon H^1(G_S,\ZM_p(r))^\chi/p^n H^1(G_S,\ZM_p(r))^\chi  \longrightarrow \left ( \bigoplus_{g\in G}  \FM_{\la^g}^\times/p^n (r-1)\right )^\chi. $$
Pour $n$ fixé on note $\LM'_n$ l'ensemble des premiers rationnels $\ell\equiv 1 [dp^n]$.
Le principal résultat original de cet article est le théorème d'indice suivant~:
\begin{theo}\label{main}
$$ \# \left (\frac {H^1(G_S,\ZM_p(r))^\chi} {\langle c^F(r)\rangle^\chi} \right )
=\max_{n\in\NM} \left ( \min_{\ell \in \LM'_n} \left (\# \frac {\left ( \bigoplus_{g\in G}  \FM_{\la^g}^\times/p^n (r-1)\right )^\chi} { \plnc (\langle c^F(r)\rangle^\chi)} \right )\right )$$
\end{theo}
Dans la section \ref{calcexp} on décrira concrètement le module $\plnc (\langle c^F(r)\rangle^\chi)$, ce qui permettra de reformuler le théorème \ref{main} de façon plus explicite.

{\it Remarque :} Un énoncé analogue au théorème \ref{main} avec des $\chi$-quotients à la place de $\chi$-parties est certainement vrai. Néanmoins, l'approche présente, qui s'appuie sur des morphismes injectifs, ne s'applique pas; et c'est l'une des raisons de choisir les $\chi$-parties pour cet article. Pour calculer les ordres $\chi$-quotients en jeu ici il faudrait développer une méthode de démonstration complètement différente.  

La suite de cette section est consacrée à la démonstration de ce théorème \ref{main}, qui passe par quelques énoncés intermédiaires intéressants en eux-mêmes. On commence par introduire des abréviations commodes.
\begin{defi} On rappelle que $F$ est abélien totalement réel de conducteur $p^a d$ avec $p\nmid d$.
\begin{enumerate}
\item Par le théorème \ref{inj}, si $n>a$ alors il existe une infinité de premiers $\ell \in \LM'_n$ tels que
le morphisme $\plnc$  soit un isomorphisme. On notera $\lnis$ cet ensemble de premiers rationnels.
\item On note $\cfrc$ le $\ZM_p[\chi]$-module libre $\cfrc=\langle c^F(r)\rangle^\chi$.
\item On note $\hsrc$ le $\ZM_p[\chi]$-module libre $\hsrc=H^1(G_S,\ZM_p(r))^\chi$.
\end{enumerate}
\end{defi}
Avec ces notations on a la relation d'indice~:
\begin{lem}\label{inegalite indices} Soit $n\in \NM$ et $\ell\in \LM'_n$. Alors
$$\# \frac{\hsrc/p^n\hsrc}{\cfrc/ (\cfrc\cap p^n \hsrc) }\ \leq \
\#\frac{\left ( \bigoplus_{g\in G}  \FM_{\la^g}^\times/p^n(r-1) \right )^\chi}  {\plnc(\cfrc /(\cfrc \cap p^n \hsrc))}.
$$
\end{lem}
\begin{proof} On avait déjà fait remarquer que les modules finis $H_{S,r}^\chi/p^n$ et  $ ( \bigoplus_{g\in G}  \FM_{\la^g}^\times/p^n(r-1)  )^\chi$ sont tous deux isomorphes à $\ZM_p[\chi]/p^n$ donc en particulier ils ont même ordre. Le lemme se ramène donc à l'inégalité évidente
$$\# \plnc(\cfrc /(\cfrc \cap p^n \hsrc)) \leq \# \cfrc /(\cfrc \cap p^n \hsrc).$$
\end{proof}
L'inégalité du lemme \ref{inegalite indices} conduit à l'égalité avec le $\displaystyle\min_{\ell\in \LM'_n}$ du théorème \ref{main}. Il suffit juste d'observer que ce min est actuellement atteint. C'est l'objet du second lemme intermédiaire qui est une conséquence immédiate du théorème \ref{inj}~:
\begin{lem}\label{indgen}
Soit $n\in \NM$.
\begin{enumerate}
\item Si $\ell\in \lnis$ alors l'isomorphisme $\plnc$ induit un isomorphisme~:
$$\frac{\hsrc/p^n\hsrc}{\cfrc/ (\cfrc\cap p^n \hsrc) }\ \overset \sim \longrightarrow \
\frac{\left ( \bigoplus_{g\in G}  \FM_{\la^g}^\times/p^n(r-1) \right )^\chi}  {\plnc(\cfrc /(\cfrc \cap p^n \hsrc))}.
$$
\item\label{indgenord} Pour tout $\ell\in \lnis$ l'inégalité d'ordre du lemme \ref{inegalite indices} est une égalité.
\item Si en outre $n$ est tel que $p^n \hsrc \subset \cfrc$ et $\ell \in \lnis$ alors $\plnc$ induit un isomorphisme~: $$\frac \hsrc \cfrc
\overset \sim \longrightarrow \
\frac{\left ( \bigoplus_{g\in G}  \FM_{\la^g}^\times/p^n(r-1) \right )^\chi}  {\plnc(\cfrc)}.
$$
\end{enumerate}
\end{lem}
\begin{proof}
Par définition même de $\lnis$ si $\ell\in \lnis$ alors l'homomorphisme $$\plnc \colon \hsrc/p^n \longrightarrow  \left ( \bigoplus_{g\in G}  \FM_{\la^g}^\times/p^n(r-1) \right )^\chi$$ est un isomorphisme. On déduit l'isomorphisme du point 1 en passant au quotient à gauche par $\cfrc/(\cfrc\cap \hsrc)$ et à droite par son image. Le point 2 est une conséquence directe du point 1. Pour démontrer le point 3, à partir des égalités $\cfrc\cap p^n\hsrc = p^n\hsrc$ et
$\cfrc+p^n\hsrc =\cfrc$, il suffit d'observer les isomorphismes évidents~:
$$ \frac \hsrc \cfrc \cong \frac \hsrc {\cfrc + p^n\hsrc} \cong \frac {\hsrc/p^n\hsrc} {\cfrc/p^n \hsrc} \cong  \frac{\hsrc/p^n\hsrc}{\cfrc/ (\cfrc\cap p^n \hsrc) } .
$$
\end{proof}
On peut énoncer le théorème ci-dessous qui permet de calculer résiduellement l'approximation modulo $p^n$ du $\chi$-indice du théorème \ref{main}~:
\begin{theo}\label{modpn}
Pour tout $n\in \NM$ on a l'égalité
$$\# \frac {\hsrc /p^n \hsrc} {\cfrc /(\cfrc \cap p^n\hsrc)} \ =\  \min_{\ell \in \LM'_n} \left ( \# \frac {\left ( \bigoplus_{g\in G}  \FM_{\la^g}^\times/p^n (r-1)\right )^\chi} { \plnc (\langle c^F(r)\rangle^\chi)} \right )
$$
\end{theo}
\begin{proof} Il suffit de combiner le lemme \ref{inegalite indices} avec le point \ref{indgenord} du lemme \ref{indgen}.\end{proof}
On démontre maintenant le théorème \ref{main}.
\begin{proof}
Comme $\hsrc$ est libre de rang $1$ sur $\ZM_p[\chi]$ qui est un anneau local, on a forcément~:
\begin{itemize}
\item soit $\cfrc \subset p^n\hsrc$, ce qui a lieu pour $n=0$ et devient faux dès que $n$ est suffisamment grand.
\item soit $p^n\hsrc \subset \cfrc$, ce qui a lieu dès que $n\geq N$ pour un certain $N$.
\end{itemize}
Cela étant, tant que $n<N$ l'approximation modulo $p^n$ ne donne aucune information puisque~:
$$\frac {\hsrc /p^n \hsrc} {\cfrc /(\cfrc \cap p^n\hsrc)} \cong \hsrc /p^n \hsrc \cong \ZM_p[\chi]/p^n. $$
Par contre dès que $n\geq N$ alors l'approximation modulo $p^n$ est exacte car on a~:
$$\frac {\hsrc /p^n \hsrc} {\cfrc /(\cfrc \cap p^n\hsrc)} \cong \frac {\hsrc /p^n \hsrc} {\cfrc /p^n \hsrc} \cong
\frac {\hsrc} {\cfrc} $$
De plus dès que $n>N$ on aura
$$\# \frac {\hsrc /p^n \hsrc} {\cfrc /(\cfrc \cap p^n\hsrc)}\ <\ \# \ZM_p[\chi]/p^n.$$
Cela explique le $\max_{n\in\NM}$ dans la formule du théorème \ref{main} et conclut sa démonstration compte-tenu du théorème \ref{modpn}.
\end{proof}
Une autre façon, peut-être plus propice aux approches numériques, de formuler le théorème \ref{main} consiste à se restreindre aux premiers rationnels $\ell\in \LM'_1$ et tel que pour au moins un $n\leq v_p(\ell-1)$, le $\chi$-indice résiduel $\displaystyle  \# \frac {\left ( \bigoplus_{g\in G}  \FM_{\la^g}^\times/p^n (r-1)\right )^\chi} { \plnc (\langle c^F(r)\rangle^\chi)} $,  soit strictement inférieur à l'indice maximal $\# \ZM_p[\chi]/p^n$. C'est une façon de comprendre l'introduction d'indices étoilés par exemple p.261 de \cite{Kur99}. De ce point de vue, on peut éviter le $\max_n$ du théorème \ref{main}. Pour cela on définit d'abord des notations supplémentaires~:
\begin{defi} On rappelle que $\LM'_1$ désigne l'ensemble des premiers rationnels $\ell\equiv 1[dp]$.
\begin{enumerate}
\item Pour tout $\ell\in \LM'_1$ et tout $n\leq v_p(\ell-1)$ on note $\mathrm{ire}_{r,\chi,\ell,n}$ (pour indice résiduel) la quantité~: $$\displaystyle {\mathrm {ire}}_{r,\chi,\ell,n} := \# \frac {\left ( \bigoplus_{g\in G}  \FM_{\la^g}^\times/p^n (r-1)\right )^\chi} { \plnc (\langle c^F(r)\rangle^\chi)} .$$
\item Pour tout $n$ on note $\mathrm{ima}_{\chi,n}$ (pour indice maximal) la quantité~: $$\displaystyle \mathrm{ima}_{\chi,n} := \# \frac {\ZM_p[\chi]} {p^n \ZM_p[\chi]}. $$
\end{enumerate}
\end{defi}
Avec ces notations, pour remplacer le maximin du théorème \ref{main} par un simple $\min$, il suffit d'écarter les cas d'égalité entre indice résiduel et indice maximal. Cela conduit au corollaire~:
\begin{cor}\label{pgcd}
$$\#\frac{H^1(G_S,\ZM_p(r))^\chi}{\langle c^F(r)\rangle^\chi}= \min_{\ell\in \LM'_1, n\leq v_p(\ell-1)} \left \{\mathrm{ire}_{r,\chi,\ell,n}; \  \mathrm {ire}_{r,\chi,\ell,n} < \mathrm{ima}_{\chi,n} \right \} .$$
\end{cor}

\section{Image explicite des $\chi$-éléments cyclotomiques}\label{calcexp}

À la vue du corollaire \ref{pgcd} il est important de rendre explicite l'indice résiduel
$$\displaystyle {\mathrm {ire}}_{r,\chi,\ell,n} =  \# \frac {\left ( \bigoplus_{g\in G}  \FM_{\la^g}^\times/p^n (r-1)\right )^\chi} { \plnc (\langle c^F(r)\rangle^\chi)}.$$
Pour ce faire on donne d'abord une formule pour un $\ZM_p[G]$-générateur  $c^F(r,\chi)$ de la $\chi$-partie $\langle c^F(r)\rangle^\chi\subset \langle c^F(r)\rangle$ et pour un générateur de $\left ( \bigoplus_{g\in G}  \FM_{\la^g}^\times/p^n (r-1)\right )^\chi$, puis on donne une formule explicite pour $\pln(c^F(r,\chi))$. Par les propriétés fonctorielles des $\chi$-parties on a l'égalité $$\plnc(\langle c^F(r)\rangle^\chi)=\langle \pln(c^F(r,\chi)) \rangle_{\ZM_p[G]}.$$
L'indice résiduel s'obtient alors en comparant $\pln(c^F(r,\chi))$ avec le générateur de $\left ( \bigoplus_{g\in G}  \FM_{\la^g}^\times/p^n (r-1)\right )^\chi$.
\begin{lem}\label{chipargen} Soit $G$ un groupe abélien et $\chi$ un caractère sur $G$. Soit $$S_\chi=\sum_{g\in \Ker \chi} g \in \ZM_p[G].$$
\begin{enumerate}
\item Si $p$ ne divise pas l'ordre $o(\chi)$ de $\chi$, alors $\ZM_p[G]^\chi\subset \ZM_p[G]$ est l'idéal principal engendré par $$T_\chi:=S_\chi \frac 1 {o(\chi)} \sum_{g\in G/\ker \chi} Tr_{\ZM_p[\chi]/\ZM_p} (\chi(g)) g^{-1} .$$
\item Si $p$ divise $o(\chi)$, soit $h\in G$ tel que $\chi(h)=\ze_p$ et soit
$\Delta\subset G/\ker \chi$ le sous-groupe des éléments d'ordre premier à $p$ de
$G/\ker \chi$. Alors $\ZM_p[G]^\chi\subset \ZM_p[G]$ est l'idéal principal
engendré par $$T_\chi:=S_\chi (1-h) \frac 1 {o(\Delta)} \sum_{\de \in \Delta} 
Tr_{\ZM_p[\chi(\Delta)]/\ZM_p} (\chi(\de)) \de^{-1}  $$
\item La $\chi$-partie $(\ZM/p^n [G])^\chi$ est le sous-module engendré par l'image mod $p^n$ du générateur  $T_\chi$ défini en 1 ou 2 ci-dessous suivant que $p$ divise ou pas $o(\chi)$.
\end{enumerate}
\end{lem}
\begin{proof} On remarque pour commencer que les générateurs $T_\chi$ sont définis indépendamment des choix de $h\in G$ tel que $\chi(h)=\ze_p$ et plus généralement des choix des relèvement des éléments de $G/\ker \chi$ dans $G$ grâce au facteur $S_\chi$ qui neutralise ces choix. Soit $\chit$ le caractère sur $G/\ker \chi$ déduit de $\chi$ par factorisation. Il est bien connu que pour tout $\ZM_p[G]$-module $M$ on a $M^\chi=(M^{\ker \chi})^\chit$. D'autre part $\ZM_p[G]^{\ker \chi} = S_\chi \ZM_p[G]$ et pour démontrer 1 et 2 on peut donc supposer sans perte que $\chi$ est fidèle sur $G$ (cyclique).
Maintenant pour voir 1 il suffit de reconnaître que l'élément $E_\chi= 1/o(G) \sum_{g\in G} Tr_{\ZM_p[\chi]/\ZM_p} (\chi(g)) g^{-1}$ est l'idempotent $\QM_p$-rationnel usuel qui définit à la fois la $\chi$-partie et le $\chi$-quotient de $\ZM_p[G]$ pour tous les caractères conjugués de $\chi$ sur $\QM_p$. Le point 1 est démontré.

Pour démontrer 2 aussi on suppose $\chi$ fidèle sur $G$ cyclique. On peut écrire $G= \Delta \times P$ où $P$ est un $p$-groupe cyclique et décomposer $\chi=\chi_1 \chi_2$ où $\chi_1$ est la restriction de $\chi$ à $\Delta$ et $\chi_2$ la restriction de $\chi$ à $P$. Pour tout $\ZM_p[G]$-module $M$ on a les égalités classiques $M^\chi=(M^{\chi_1})^{\chi_2}=(E_{\chi_1} M)^{\chi_2}$ et on se ramène donc à démontrer 2 dans le cas où $G=P$ est un $p$-groupe cyclique et $\chi$ est fidèle sur $G$. Dans ce cas l'unique $h\in G$ tel que $\chi(h)=\ze_p$ est un générateur de l'unique sous-groupe $C\subset G$ d'ordre $p$. Et par le lemme II.2 de \cite{So90} on a pour tout $\ZM_p[G]$-module $M$ l'égalité $M^\chi=\ker N_C$ où $N_C\colon M\longrightarrow M$ est l'application norme suivant $C$ définie par $N_C(m)= \sum_{g\in H} g m$.
Dans le cas particulier où $M=\ZM_p[G]$ alors $M$ est $C$-cohomologiquement trivial et on a $\ker N_C = (1-h) \ZM_p[G]$, ce qui conclut la preuve de 2.

Pour voir 3 on définit pour tout $\ZM_p[G]$-module $M$ le module "étendu-tordu" $\RC_{\chi,M}$ comme le $\ZM_p[\chi]$-module $\RC_{\chi,M}=M\otimes \ZM_p[\chi]$ muni de l'action de $G$ tordue $g* x = \chi(g^{-1}) g x$. Une minute de réflexion permet de s'assurer de l'isomorphie $M^\chi\cong (\RC_{\chi,M})^G$. Ceci posé, en partant de la suite $$\xymatrix {0 \ar[r] & p^n \ZM_p[G] \ar[r] & \ZM_p[G] \ar[r] & \ZM/p^n[G] \ar[r] & 0},$$ on obtient en passant aux $\chi$-parties et en utilisant 1 ou 2 suivant que $p$ divise ou pas $o(\chi)$ la suite~: $$\xymatrix {0 \ar[r] & p^n T_\chi \ZM_p[G] \ar[r] & T_\chi \ZM_p[G] \ar[r] & (\ZM/p^n[G])^\chi \ar[r] & H^1(G,\RC_{\chi,p^n\ZM_p[G]})}.$$
Mais pour ce qui concerne le module libre $p^n \ZM_p[G]$ son étendu-tordu $ \RC_{\chi,p^n\ZM_p[G]}$ est encore isomorphe à $\ZM_p[\chi][G]$ qui est $G$-cohomologiquement trivial, d'où l'isomorphie requise $ (T_\chi \ZM_p[G])/p^n \cong (\ZM/p^n[G])^\chi$.
\end{proof}
Par le théorème \ref{freespe} le module $\langle c^F(r)\rangle$ est $\ZM_p[G]$-libre et l'on déduit immédiatement du lemme \ref{chipargen} l'égalité $\langle c^F(r)\rangle^\chi= \langle T_\chi c^F(r)\rangle$. D'où la proposition~:
\begin{prop} Soit $\chi$ un caractère sur $G$, soit $\De\subset G/\ker \chi$ le sous-groupe des éléments d'ordre premier à $p$ de $G/\ker \chi$ et si $p\mid o(\chi)$ soit $h\in G$ tel que $\chi(h)=\ze_p$. La $\chi$-partie $\langle c^F(r)\rangle^\chi$ est le sous-$\ZM_p[G]$-module monogène engendré par
$$ T_\chi c^F(r) =\left \{ \begin{array}{l}
                             \displaystyle \frac 1 {o(\chi)} \sum_{g\in G} Tr_{\ZM_p[\chi]/\ZM_p} (\chi(g^{-1})) g c^F(r) \ \text{si} \ p\nmid o(\chi) \\
                              \displaystyle \frac {1-h} {o(\De)} \sum_{g\in G, p\nmid o(\chi(g))} Tr_{\ZM_p[\chi]/\ZM_p} (\chi(g^{-1})) g c^F(r) \ \text{si}\ p\mid o(\chi)  \\
                            \end{array}
                          \right.$$
\end{prop}
\begin{proof} Les égalités sont de simples ré-écritures de l'élément $T_\chi$ appliqué à $c^F(r)$.
\end{proof}

\begin{prop}\label{explipln} Soit $\ell \equiv 1 [dp^n]$. On fixe $\xi_\ell\in \ZM$ une racine primitive modulo $\ell$ telle que $\ze_{dp^n} \equiv \xi_\ell^{(\ell-1)/dp^n} [\LC(\la)]$, où $\la$ et $\LC(\la)$ sont les idéaux divisant $\ell$ fixé pour la construction de $\phi_{\ell,n}$ en début de section \ref{redl}. Partant de l'isomorphisme canonique $\Gal(\QM(\ze_{dp^n})/\QM) \cong (\ZM/dp^n)^\times$, on fixe $I_{F,n}$ un système de représentants dans $\ZM$ de $\Gal(\QM(\ze_{dp^n})/F)$ et, pour tout $g\in G$, un représentant $a_g\in \ZM$ de $g^{-1}$. On a~:
$$\pln (c^F(r))= \left ( \prod_{i\in I_{F,n}} (1-\xi_\ell^{a_g i(\ell-1)/(dp^n)})^{i^{r-1}}\otimes t(r-1)_n\right )_{g\in G}.$$
\end{prop}
\begin{proof} On a choisit $\xi_\ell$ pour avoir $\displaystyle \ze_{dp^n} \equiv \xi_\ell^{(\ell-1)/dp^n}[\LC(\la)]$, et donc aussi $\displaystyle \ze_{dp^n}^{a_g} \equiv \xi_\ell^{a_g (\ell-1)/dp^n}[\LC(\la)]$. Pour tout $g\in G$ soit  $\si_g\in \Gal(\QM(\ze_{dp^n})/\QM)$ le relèvement de $g$ tel que $\displaystyle \ze_{dp^n}=\ze_{dp^n}^{\si_g a_g}$ alors on a
$\displaystyle \ze_{dp^n} =\ze_{dp^n}^{a_g \si_g} \equiv \xi_\ell^{a_g (\ell-1)/dp^n}[\LC(\la)^{\si_g}]$, parce que $\si_g$ agit trivialement sur $\xi_\ell$. On calcule maintenant $\pln(c^F(r))$~:
$$\begin{aligned}\pln ( c^F(r))&=    \left ( res_{\la^g} ( \varphi_n (c^F(r)))\right )_{g\in G}\ \text{par définition de}\ \pln,\\
&=   \left ( red_{\la^g}(r-1) ( \varphi_n ( c^F(r)))\right )_{g\in G}\ \text{voir diagramme (\ref{kumres}})\\
&=   \left ( red_{\la^g}(r-1) \left ( \sum_{\si}  \kappa^{r-1}(\si) \si (1-\ze_{dp^n}) \otimes t(r-1)_n \right)\right )_{g\in G}\\ &  \text{par la formule (\ref{phi(c)}) et où}\ \si\ \text{parcourt}\ \Gal(\QM(\ze_{dp^n})/F) \\
&=  \left ( red_{\la^g}(r-1) ( (\prod_{i\in I_{F,n}} (1-\ze_{dp^n}^i)^{i^{r-1}}) \otimes t(r-1)_n) \right )_{g\in G}\\
&=  \left ( \prod_{i\in I_{F,n}} (1-\xi_\ell^{a_g i(\ell-1)/(dp^n)})^{i^{r-1}}\otimes t(r-1)_n\right )_{g\in G}.\\
\end{aligned}$$
En effet la congruence $\displaystyle \ze_{dp^n}\equiv \xi_\ell^{a_g (\ell-1)/dp^n}[\LC(\la)^{\si_g}]$ dans $\QM(\ze_{dp^n})$ remarquée en début de preuve donne dans $F^\times/p^n (r-1)$~:
$$red_{\la^g}(r-1) (\prod_{i\in I_{F,n}} (1-\ze_{dp^n}^i)^{i^{r-1}}\otimes t(r-1)_n) = \prod_{i\in I_{F,n}} (1-\xi_\ell^{a_g i(\ell-1)/(dp^n)})^{i^{r-1}} \otimes t(r-1)_n). $$
 Il faut aussi noter que la formule finale ne dépend plus des choix de $a_g$ ou de $\si_g$ puisque $\Gal(\QM(\ze_{dp^n}/F)$ agit trivialement de part et d'autre.
\end{proof}
On a maintenant tous les éléments pour donner une formule explicite calculant les quantités $\mathrm{ire}_{r,\chi,\ell,n}$ qui interviennent dans le corollaire \ref{pgcd}.
\begin{theo}\label{formexpire} Soit $\ell\equiv 1[dp]$ et soit $n\leq v_p(\ell-1)$. Pour tout groupe abélien $M$ fini on note $M[p^n]$ les éléments de $M$ annulés par $p^n$. Soit $\eta_\ell$ une racine primitive modulo $\ell$. On a
\begin{equation}\label{formexpireeq}\mathrm{ire}_{r,\chi,\ell,n}=\# \frac {T_\chi \left ( \bigoplus_{g\in G}  \FM_{\la^g}^\times [p^n]\right )}
{ \langle T_\chi  \left(\prod_{i\in I_{F,n}}
(1-\eta_\ell^{a_g i(\ell-1)/(dp^n)} )^{i^{r-1}\frac {\ell-1} {p^n}}\right)_{g\in G}\rangle_{\ZM_p[\chi]} }
\end{equation}
\end{theo}
\begin{proof}
Pour calculer l'indice $$\# \frac {T_\chi \left ( \bigoplus_{g\in G}  \FM_{\la^g}^\times/p^n (r-1)\right )}
{ \langle T_\chi \left ( \prod_{i\in I_{F,n}} (1-\xi_\ell^{a_g i(\ell-1)/(dp^n)})^{i^{r-1}}\otimes t(r-1)_n\right )_{g\in G}\rangle } $$ on commence par enlever les $r-1$ tordus à la Tate au numérateur et au dénominateur ce qui ne change pas l'indice. 
Ensuite si on prend $\eta_\ell=\xi_\ell$ on trouve d'après la proposition \ref{explipln} et le lemme \ref{chipargen} l'égalité
$$\mathrm{ire}_{r,\chi,\ell,n}=\# \frac {T_\chi \left ( \bigoplus_{g\in G}  \FM_{\la^g}^\times /p^n\right )}
{ \langle T_\chi  \left(\prod_{i\in I_{F,n}}
(1-\eta_\ell^{a_g i(\ell-1)/(dp^n)} )^{i^{r-1}}\right)_{g\in G}\rangle_{\ZM_p[\chi]} }. $$
Clairement cette égalité reste valable avec tout choix de racine primitive modulo $\ell$ autre que $\xi_\ell$ parce que l'indice de gauche dans l'égalité (\ref{formexpireeq}) ne dépend pas de ce choix.
L'exponentiation par $(\ell-1)/p^n$ réalise un isomorphisme entre $\FM_\la^\times/p^n$ et $\FM_\la^\times [p^n]$, ce qui conclut la preuve du théorème \ref {formexpire}.
\end{proof}
La formule finale (\ref{formexpireeq}) s'entend comme l'indice de deux $\ZM_p[\chi]$-modules monogènes. Dans le cas particulier où l'ordre de $\chi$ divise $p-1$ alors $\ZM_p[\chi]$ s'identifie à $\ZM_p$. Il est alors possible, comme c'est fait dans \cite{Kur99} p. 262, de définir des $c_{l,r,p^n}^\chi\in \FM_l^\times$ et non dans un produit de $\FM_\la^\times$; puis de calculer les $\mathrm{ire}_{r,\chi,l,n}$ comme des indices dans un seul $\FM_l^\times$. Dans notre cadre plus général on peut avoir besoin de $[\QM_p(\chi):\QM_p]$ paramètres $p$-adiques et on est contraint de travailler dans $ \bigoplus_{g\in G}  \FM_{\la^g}^\times $. Cela augmente la complexité mais ne diminue pas l'effectivité de ces calculs d'indices.

\section{Quelques r\'esultats num\'eriques}\label{numerique}

Voici, suivant une des suggestions du referee, quelques exemples numériques obtenus avec le logiciel SAGE. Le code SAGE utilisé et d'autres données numériques sont disponible sur la page http://www-irma.u-strasbg.fr/$\sim$beliaeva/codes$.$html. Pour  simplifier l'algorithme, tous les cas présentés ici concernent des corps 
de nombre de degr\'e $p$. Pour ces corps il n'y a qu'un caractère non trivial à $\QM_p$-conjugaison près, d'o\`u l'absence de variation en fonction du choix des caract\`eres. On se restreint aussi à des corps de conducteur $d$ premier avec  $d\equiv 1 [p]$, pour que $\QM(\zeta_d)$ contienne un unique sous-corps de degr\'e $p$ sur $\QM$. On note $F_{d,p}$ ce sous-corps de degré $p$ de $\QM(\ze_d)$. D'après le corollaires \ref{pgcd}, pour $F=F_{d,p}$ on a~:
$$\#\frac{H^1(G_S,\ZM_p(r))^\chi}{\langle c^{F_{d,p}}(r)\rangle^\chi}= \min_{\ell\in \LM'_1, n\leq v_p(\ell-1)} \left \{\mathrm{ire}_{r,\chi,\ell,n}; \  \mathrm {ire}_{r,\chi,\ell,n} < \mathrm{ima}_{\chi,n} \right \} .$$
On a donc calculé les $\mathrm{ire}_{r,\chi,\ell,n}$ pour $p=3,5,7,11$ et  $\ell<10^6$. Cela donne une majoration de l'indice $\#\frac{H^1(G_S,\ZM_p(r))^\chi}{\langle c^{F_{d,p}}(r)\rangle^\chi}$, mais au vu des résultats numériques on peut s'attendre à ce que le $\mathrm{ire}_{r,\chi,\ell,n}$ retenu soit égal à cet indice. En tout cas cela se produit chaque fois qu'un $\mathrm{ire}_{r,\chi,\ell,n}$ est égal à $p$.

En effet l'indice calculé $\#\frac{H^1(G_S,\ZM_p(r))^\chi}{\langle c^F(r)\rangle^\chi}$ n'est jamais trivial à cause de la différence entre 
$C^F(r)$ et $\langle c^F(r)\rangle$. Concrètement, avec la relation (\ref{reltor}) entre $c^F(r)$ et l'autre générateur $c_1^F(r)$ du module $C^F(r)$, on obtient 
que le quotient $C^F(r)/\langle c^F(r)\rangle$ est isomorphe à $\ZM_p/(1-d^{r-1})\ZM_p$ avec action triviale de $G$. Et comme $p\mid(1-d^{r-1})$, on trouve pour tout  caractère $\chi$ non trivial sur $G$ l'isomorphisme $(C^F(r)/\langle c^F(r)\rangle)^\chi \simeq \ZM_p/p\ZM_p$. En conséquence si on part de la suite 
\(\xymatrix {0\ar[r] &\frac {C^F(r)} {\langle c^F(r) \rangle} \ar[r]  & \frac{H^1(G_S,\ZM_p(r))}{\langle c^F(r)\rangle} \ar[r] & \frac{H^1(G_S,\ZM_p(r))} {C^F(r)} \ar[r] & 0},\) on obtient en passant aux $\chi$-parties la suite 
\(\xymatrix {0\ar[r] &\ZM_p/p\ZM_p  \ar[r]  & \frac{H^1(G_S,\ZM_p(r))^\chi}{\langle c^F(r)\rangle^\chi} \ar[r] & \left (\frac{H^1(G_S,\ZM_p(r))} {C^F(r)}\right)^\chi \ar[r] & \cdots}.\) Ainsi dès que l'indice calculé ${\mathrm {ire}}_{r,\chi,\ell,n}$ est exactement $p$ on sait que c'est aussi l'indice $\#\frac{H^1(G_S,\ZM_p(r))^\chi}{\langle c^F(r)\rangle^\chi}$.

Pour chacun de ces corps on a testé les $r$-tordus à la Tate pour des $r$ impairs allant de 3 \`a 21. 
Dans certains cas il n'y a pas de variation en fonction de $r$. Le m\^eme indice appara\^it pour la premi\`ere fois pour le m\^eme couple $(\ell,n)$.

Voici quelques exemples où le résultat des calculs effectués ne dépend pas de $r$; dans tous ces cas on remarque que l'indice minimal appara\^it pour le plus petit $\ell$ dans $\LM'_1$~: 
\begin{enumerate}
	\item $p=3$, $d=7$ ($F=F_{3,7}=\QM(\cos(\frac{2\pi}7)$): l'indice est $3$, il appara\^it pour la premi\`ere fois pour $\ell=43,\ n=1$ ;
	\item $p=3$, $d=13$ ($F=L_{3,13}=\QM( \cos(\frac{2\pi}{13})+\cos(\frac{10\pi}{13}))$): l'indice est $3$, il appara\^it pour la premi\`ere fois pour $\ell=79,\ n=1$;
	\item $p=5$, $d=11$ ($F=F_{5,11}=\QM(\cos(\frac{2\pi}{11}))$): l'indice minimal trouv\'e est $5^3$, il appara\^it pour la premi\`ere fois pour $\ell=331,\ n=1$;
	\item $p=5$, $d=41$ ($F=F_{5,41}=\QM(\cos(\frac{2\pi}{41})+\cos(\frac{3\pi}{41})+ \cos(\frac{9\pi}{41})+ \cos(\frac{14\pi}{41})$): l'indice minimal trouv\'e est $5^3$, il appara\^it pour la premi\`ere fois pour $\ell=821,\ n=1$;
	\item $p=7$, $d=29$ ($F=F_{7,29}$): l'indice minimal trouv\'e est $7^5$, il appara\^it pour la premi\`ere fois pour $\ell=2437,\ n=1$;
	\item $p=7$, $d=43$ ($F=F_{7,43}$): l'indice minimal trouv\'e est $7^5$, il appara\^it pour la premi\`ere fois pour $\ell=3011,\ n=1$;
	\item $p=11$, $d=23$ ($F=F_{11,23}$): l'indice minimal trouv\'e est $11^9$, il appara\^it pour la premi\`ere fois pour $\ell=1013,\ n=1$.
\end{enumerate}
Voici d'autres exemples où l'indice d\'epend de $r$ :
\begin{table}[h]
	\centering
		\begin{tabular}[t]{|l|c|c|c|c|}
		\hline
			$p$&$d$&$r$&indice& premi\`ere apparition\\
		\hline
			3&19& 3,5,9,11,15,17,21 & $3^3$& $\ell=2053$, $n=2$\\
		\cline{3-5}
		   & & 7,13& $3^5$ & $\ell=2053$, $n=3$\\
		\cline{3-5}
		   & & 19& $3^7$& $\ell=3079$, $n=4$\\
		\cline{2-5}
				&37& 3,5,9,11,15,17,21 & $3^5$ & $\ell=1999$, $n=3$\\
		\cline{3-5}
			& &7,13& $3^7$& $\ell=41959$, $n=4$\\
		\cline{3-5}
			& &19& $3^9$ &$ \ell= 161839$, $n=5$  \\
		\hline
		  5&31& 3,7,11,15,19 &$5^3$& $\ell=311$, $n=1$\\
		\cline{3-5}
		  & & 5,9,13,17& $5^7$& $\ell=4651$, $n=2$\\
		\cline{3-5}
		  & & 21& $5^{11}$& $\ell=23251$, $n=3$\\
		\hline
	\end{tabular}
\end{table}

Dans les exemples trait\'es pour  $p=7$ et $p=11$ on n'a pas de variation en fonction de $r$.

\backmatter


\begin{thebibliography}{10}
\bibitem{Bei83} {\scshape A.~A. Be{\u\i}linson} -- {\og Higher regulators of modular
  curves\fg}, Applications of algebraic {$K$}-theory to algebraic geometry and
  number theory, {P}art {I}, {II} ({B}oulder, {C}olo., 1983), Contemp. Math.,
  vol.~55, Amer. Math. Soc., Providence, RI, 1986, p.~1--34.

\bibitem{Tan} {\scshape T.~Beliaeva} -- {\og Unit\'{e}s semi-locales modulo sommes de {G}au\ss\ en th\'eorie d'{I}wasawa\fg}, \emph{Th\`{e}se de l'universit\'{e} de Franche-Comt\'{e} Besan\c{c}on} (2004).

\bibitem{CJM} {\scshape J.-R. Belliard} -- {\og Global units modulo circular units: descent without {I}wasawa's main conjecture\fg}, \emph{Canad. J. Math.} \textbf{61} (2009), no.~3, p.~518--533.

\bibitem{BenN02} {\scshape D.~Benois {\normalfont \smfandname} T.~{N}guy{\cftil{e}}n-{Q}uang {D}{\cftil{o}}} -- {\og Les nombres de {T}amagawa locaux et la conjecture de {B}loch et {K}ato pour les motifs {$\Bbb Q(m)$} sur un corps ab\'elien\fg}, \emph{Ann. Sci. \'Ecole Norm. Sup. (4)} \textbf{35} (2002), no.~5, p.~641--672.

\bibitem{BK90} {\scshape S.~Bloch {\normalfont \smfandname} K.~Kato} -- {\og {$L$}-functions and {T}amagawa numbers of motives\fg}, The Grothendieck Festschrift, Vol.\ I, Progr. Math., vol.~86, Birkh\"auser Boston, Boston, MA, 1990, p.~333--400.

\bibitem{BG03} {\scshape D.~Burns {\normalfont \smfandname} C.~Greither} -- {\og On the equivariant {T}amagawa number conjecture for {T}ate motives\fg}, \emph{Invent. Math.} \textbf{153} (2003), no.~2, p.~303--359.

\bibitem{Col79} {\scshape R.~F. Coleman} -- {\og Division values in local fields\fg}, \emph{Invent. Math.} \textbf{53} (1979), no.~2, p.~91--116.


\bibitem{D89} {\scshape P.~Deligne} -- {\og Le groupe fondamental de la droite projective moins trois points\fg}, Galois groups over ${\bf Q}$ (Berkeley, CA, 1987), Math. Sci. Res. Inst. Publ., vol.~16, Springer, New York, 1989, p.~79--297.

\bibitem{Ga01} {\scshape W.~Gajda} -- {\og On cyclotomic numbers and the reduction map for the {$K$}-theory of the integers\fg}, \emph{$K$-Theory} \textbf{23} (2001), no.~4, p.~323--343.

\bibitem{Grei92} {\scshape C.~Greither} -- {\og Class groups of abelian fields, and the main conjecture\fg}, \emph{Ann. Inst. Fourier (Grenoble)} \textbf{42} (1992), no.~3, p.~449--499.

\bibitem{HK03} {\scshape A.~Huber {\normalfont \smfandname} G.~Kings} -- {\og Bloch-{K}ato conjecture and {M}ain {C}onjecture of {I}wasawa theory for {D}irichlet characters\fg}, \emph{Duke Math. J.} \textbf{119} (2003), no.~3, p.~393--464.

\bibitem{HuW} {\scshape A.~Huber {\normalfont \smfandname} J.~Wildeshaus} -- {\og Classical motivic polylogarithm according to {B}eilinson and {D}eligne\fg}, \emph{Doc. Math.} \textbf{3} (1998), p.~27--133 (electronic).

\bibitem{Kosur} {\scshape M.~Kolster} -- {\og {$K$}-theory and arithmetic\fg}, Contemporary developments in algebraic {$K$}-theory, ICTP Lect. Notes, XV, Abdus Salam Int. Cent. Theoret. Phys., Trieste, 2004, p.~191--258 (electronic).

\bibitem{KNF} {\scshape M.~Kolster, T.~{N}guy{\cftil{e}}n {Q}uang {D}{\cftil{o}} {\normalfont \smfandname} V.~Fleckinger} -- {\og Twisted ${S}$-units, $p$-adic class number formulas, and the {L}ichtenbaum conjectures\fg}, \emph{Duke Math. J.} \textbf{84} (1996), no.~3, p.~679--717.

\bibitem{Kur99} {\scshape M.~Kurihara} -- {\og The {I}wasawa {$\lambda$}-invariants of real abelian fields and the cyclotomic elements\fg}, \emph{Tokyo J. Math.} \textbf{22} (1999), no.~2, p.~259--277.

\bibitem{NSW} {\scshape J.~Neukirch, A.~Schmidt {\normalfont \smfandname} K.~Wingberg} -- \emph{Cohomology of number fields}, Grundlehren der Mathematischen Wissenschaften [Fundamental Principles of Mathematical Sciences], vol. 323, Springer-Verlag, Berlin, 2000.

\bibitem{Kosur75} {\scshape D.~Quillen} -- {\og Higher algebraic {$K$}-theory. {I}\fg}, Algebraic {$K$}-theory, {I}: {H}igher {$K$}-theories ({P}roc. {C}onf., {B}attelle {M}emorial {I}nst., {S}eattle, {W}ash., 1972), Springer, Berlin, 1973, p.~85--147. Lecture Notes in Math., Vol. 341.

\bibitem{Si78} {\scshape W.~Sinnott} -- {\og On the {S}tickelberger ideal and the circular units of a cyclotomic field\fg}, \emph{Ann. of Math. (2)} \textbf{108} (1978), no.~1, p.~107--134.

\bibitem{Si80} \bysame , {\og On the {S}tickelberger ideal and the circular units of an abelian field\fg}, \emph{Invent. Math.} \textbf{62} (1980), no.~2, p.~181--234.

\bibitem{So90} {\scshape D.~Solomon} -- {\og On the classgroups of imaginary abelian fields\fg}, \emph{Ann. Inst. Fourier (Grenoble)} \textbf{40} (1990), p.~467--492.

\bibitem{Sou80} {\scshape C.~Soul{\'e}} -- {\og On higher $p$-adic regulators\fg}, Algebraic $K$-theory, Evanston 1980 (Proc. Conf., Northwestern Univ., Evanston, Ill., 1980), Springer, Berlin, 1981, p.~372--401.

\bibitem{Sou87} \bysame , {\og \'{E}l\'ements cyclotomiques en {$K$}-th\'eorie\fg}, \emph{Ast\'erisque} (1987), no.~147-148, p.~225--257, 344, Journ\'ees arithm\'etiques de Besan\c con (Besan\c con, 1985).

\bibitem{T99} {\scshape T.~Tsuji} -- {\og Semi-local units modulo cyclotomic units\fg}, \emph{J. Number Theory} \textbf{78} (1999), no.~1, p.~1--26. 

\bibitem{V10} {\scshape V.~Voevodsky} -- {\og On motivic cohomology with $\ZM/l$-coefficients
\fg}, \emph{Ann. of Math. (2)} \textbf{174} (2011), no.~1, p.~401--438. 

\end{thebibliography}
\end{document}